\documentclass[12pt]{article}

\usepackage{amssymb}
\usepackage{amsfonts}
\usepackage{latexsym}
\usepackage{graphicx}

\usepackage[usenames]{color}
\usepackage{amsmath,amsthm}

\usepackage{epsfig}

\usepackage{color}
\usepackage{amsmath,amsbsy,amsfonts,amssymb,latexsym,graphics,epsfig}

\newcommand{\ignore}[1]{}

\newcommand{\beqn}{\begin{eqnarray*}}
\newcommand{\eeqn}{\end{eqnarray*}}

\newcommand{\hd}{\mathcal{H}}
\newcommand{\tl}{\mathcal{T}}

\newcommand{\ZZ}{\mathcal{Z}}

\newcommand{\eps}{\varepsilon}

\newcommand{\ee}{\mathrm{e}}

\newcommand{\Z}{\mathbb{Z}}

\newcommand{\R}{\mathbb{R}}

\newcommand{\shf}{\mbox{\footnotesize $\frac{1}{2}$}}

\newcommand{\dd}{\mathrm{d}}
\newcommand{\DD}{\mathrm{D}}

\newcommand{\Matrix}[1]{\ensuremath{\left[\begin{array}{ccccccccccccccccr} #1 \end{array}\right]}}

\newcommand{\etal}{{\em et al.}}

\newcommand{\CC}{{\mathcal C}}
\newcommand{\EE}{{\mathcal E}}

\newcommand{\GG}{{\mathcal G}}

\newcommand{\KK}{{\mathcal K}}
\newcommand{\MM}{{\mathcal M}}

\newcommand{\mS}{{\mathcal S}}
\newcommand{\XX}{{\mathcal X}}
\newcommand{\YY}{{\mathcal Y}}

\newcommand{\QQ}{{\mathcal Q}}

\newtheorem{theorem}{Theorem}[section]
\newtheorem{lemma}[theorem]{Lemma}
\newtheorem{definition}[theorem]{Definition}
\newtheorem{corollary}[theorem]{Corollary}
\newtheorem{proposition}[theorem]{Proposition}
\newtheorem{remark}[theorem]{Remark}
\newtheorem{remarks}[theorem]{Remarks}
\newtheorem{example}[theorem]{Example}

\renewcommand{\epsilon}{\varepsilon}
\newcommand{\Sone}{{\mathbb S}^1}

\newcommand{\ii}{\mbox{i}}

{\begingroup

  \begin{enumerate}}
  {\end{enumerate}\endgroup}

\title{Stable Synchronous Propagation of Periodic Signals \\
by Feedforward Networks}
\author{Ian Stewart  and David Wood\\ Mathematics Institute\\ University of Warwick \\ Coventry CV4 7AL
\\ United Kingdom}

\date{\today}

\begin{document}

\maketitle

\begin{abstract}
We analyse the dynamics of networks in which a central pattern generator
(CPG) transmits time-periodic signals along one or more feedforward chains in a
synchronous or phase-synchronous manner. Such propagating signals are common in
biology, especially in locomotion and peristalsis, and are of interest for
continuum robots. We construct such networks 
as feedforward lifts of the CPG. If the CPG dynamics is periodic, so is the lifted dynamics. 
Synchrony with the CPG manifests as a standing wave, and a regular phase 
pattern creates a travelling wave. We discuss Liapunov, asymptotic,
 and Floquet stability of
the lifted periodic orbit and introduce transverse versions of these conditions
that imply stability for signals propagating along arbitrarily long chains.
We compare these notions to a simpler condition,
transverse stability of the synchrony subspace,
which is equivalent to Floquet stability when nodes are $1$-dimensional. 
\end{abstract}

\section{Introduction}
\label{S:intro}

Many aspects of animal physiology involve the longitudinal propagation of 
rhythmic time-periodic patterns in which linear chains of neurons oscillate in synchrony 
or with specific phase relations. These two types of
behaviour can be interpreted as standing waves and travelling waves, respectively.
A common mechanism for such propagating chains involves
a network of neurons, often called a Central Pattern Generator (CPG),
which generates the basic rhythms. This lies at the start of a feedforward
network along which the CPG signals propagate. Similar waves of motion are used
to propel snake-like robots for exploration (including other planets) \cite{MWXLW17}; 
there are also numerous medical applications, see
\cite{JC19,SOWRK10,ZHX20} and references therein. This field of `continuum robots' is
advancing rapidly and the literature is huge.

Both types of application can be modelled using networks of coupled 
dynamical systems. We work in the general formalism of \cite{GS23,GST05,SGP03},
see Section~\ref{S:NAO}. 
We say that two nodes are {\em synchronous} if their waveforms
(time series) are identical. More generally, two nodes are {\em phase-synchronous} if their waveforms
(time series) are identical except for a phase shift (time translation).
These definitions are idealisations, but they open up a
powerful mathematical approach with useful implications. Real
systems can be considered as perturbations of
idealised ones, and much of the interesting structure can persist in an appropriately
approximate form. 

The main aim of this paper is to
describe a general method for constructing networks in which periodic
dynamics of a specified CPG propagates synchronously, or phase-synchronously
with a regular pattern of phase shifts, along a feedforward chain, 
tree, or any other feedforward structure. (For simplicity we often use the term `chain'
without implying linear topology.) This is achieved by 
constructing the rest of the network as a {\em feedforward lift} of the CPG.

Of course, the use of chains to propagate signals is not a new idea,
as even a cursory glance at the literature shows. Indeed,
it is arguably the simplest, most natural, and most obvious method. However, the
formal setting in which we carry out the analysis makes it possible to prove some 
general stability results and helps to unify the area. 

\subsubsection{Stability}

A key issue is to ensure that these propagating states are stable.
This term has many technically different meanings, see \cite{BS70} and
Section~\ref{S:BStab}.
More recently, chaotic dynamics has extended the diversity of meanings. 
Stability of synchronous states has been widely studied for
special models, such as the Kuramoto model \cite{K84,MBBP19}.
Other approaches and related results can be found in \cite{BPP00,PC98,PK17},
and a version for random dynamical systems is analysed in \cite{HLC13}.

We consider several notions
of stability for equilibria and periodic orbits, concentrating on the periodic case.
Roughly speaking, Liapunov stability means that a small perturbation of the initial conditions
has a small effect on the orbit; asymptotic stability means that the state converges to
the orbit after a small perturbation; and for exponential stability
the convergence has an exponential bound. For an equilibrium, exponential stability is
equivalent to linear stability; for a periodic orbit
it is equivalent to stability in the Floquet sense \cite{HKW81}, which
for brevity we call `Floquet stability'. For formal definitions and further
discussion, see Section \ref{S:BStab}.

\subsubsection{Transverse Stability}

The feedforward structure implies that if the lifted periodic orbit
is stable, for a given stability notion, then the CPG orbit must be stable.
However, this condition is not sufficient for stability of the lifted periodic orbit,
because synchrony might be destroyed by perturbations transverse to the
synchrony subspace, that is, by {\em synchrony-breaking} perturbations.
The main point of this paper is to find necessary and sufficient 
conditions for the lifted state to be stable in each of the three senses above.
This is achieved by defining associated conditions of `transverse' 
Liapunov, asymptotic, and linear/Floquet stability.

For each of the three stability notions $\mS$, we prove that the
lifted periodic state is $\mS$-stable if and only if the CPG periodic orbit is
$\mS$-stable on the CPG state space and the orbit is transversely $\mS$-stable
at every node of the chain. These results are stated and proved in
Theorem \ref{T:FFStab} for Floquet stability and in 
Theorem \ref{T:FFLSstab} for Liapunov stability. There
is also a version for asymptotic stability; we omit a statement and proof
since these are similar to, and simpler than, those for Liapunov stability.
Theorem \ref{T:TWtranseigen} generalises the Floquet stability result to signals
that propagate according to a specified phase pattern.

A related, and simpler, notion is transverse stability of the synchrony subspace. 
Intuitively, this states that the vector field is
attracting towards the synchrony subspace at every point
on that subspace. Technically, it means that at any point on the synchrony subspace
(or, more generally, in a neighbourhood of the CPG periodic orbit) all
eigenvalues of the Jacobian,
for eigenvectors transverse to the synchrony subspace, 
have negative real parts.  If
the state of the CPG is Floquet stable and node spaces are 1-dimensional,
this condition implies Floquet stability of the lifted state. For node spaces
of dimension 2 or more, this implication trivially remains valid for
equilibria, but it can fail for periodic orbits, as the celebrated Markus-Yamabe counterexample (Example \ref{ex:notFloquet}) shows. Despite this,
transverse stability of the synchrony subspace
retains some heuristic value, and can sometimes be 
given rigorous justification. It is therefore worth examining in its own right,
independently of its relation to overall stability.

An important feature of feedforward lifts is that
transverse Floquet stability is determined by dynamics associated with
individual nodes of the CPG. In consequence, our results
show that if the propagating signal is Floquet stable one step along the chain, then
it remains Floquet stable however long the chain is, or if the chain branches like a tree.
A side effect of this feature is the generic occurrence of multiple
Floquet multipliers, except for very short chains, even when the overall network
has no symmetry. This multiplicity is an advantage for all forms of transverse stability,
but may cause problems for bifurcation analysis when it occurs for a critical eigenvalue.

Similar remarks apply to transverse Liapunov stability, but now it is necessary to consider
states sufficiently close to the periodic one, not just on it. We plan to discuss this point in a
future paper.

Some of the results generalise to the propagation of signals that are
not periodic, in particular those involving Liapunov stability: see Section \ref{S:SCS}. 
For simplicity we focus on the time-periodic case, where the stronger notion of 
Floquet stability applies. In this setting we also discuss propagation of phase-related
signals, where certain nodes have the same waveform subject to regular phase shifts.
Such signals can be viewed as travelling waves.

\subsection{Biological Motivation}
\label{SS:BM}

To set the scene, we begin with two examples of propagating phase-synchronous signals
in biological systems: peristalsis in the gut and peristaltic waves in 
crawling movement in {\em Drosophila} larvae.
Further examples include
the heartbeat of the medicinal leech \cite{BP04,CP83,CNO95}, 
legged locomotion \cite{B01,B19,CS93a,CS93b,GSBC98,GSCB99,PG06,S14},
and the motion of the nematode worm
{\em Caenorhabditis elegans} \cite{BBC12,IB18,OIB21,SSSIT21}.
These networks are similar, but not identical, to feedforward lifts,
and are presented solely as motivation.

\begin{example}\em
\label{ex:persistalsis}
Peristalsis in the intestine is a travelling wave of muscular contractions
controlled by the enteric nervous system, which 
contains millions of neurons, mainly bunched into 
ganglia of two types: myenteric and submucosal.
Successive ganglia are connected together, and the large scale topology
for each type is that of a chain.
Submucosal ganglia are spaced more closely than myenteric ones.

General information is in \cite{F08,Gr03}.
Mathematical models of enteric neural motor patterns are surveyed in~\cite{CTB13},
which contains few mathematical details but a large number of references.
We also mention \cite{CBT08} on a model of intestinal segmentation
and \cite{TBB00} on a recurrent excitatory network model, both in guinea pigs.
 Figure~\ref{F:KF99fig}
shows a schematic network from \cite{KF99} with modular feedforward structure.

\begin{figure}[h!]
\centerline{%
\includegraphics[width=0.4\textwidth]{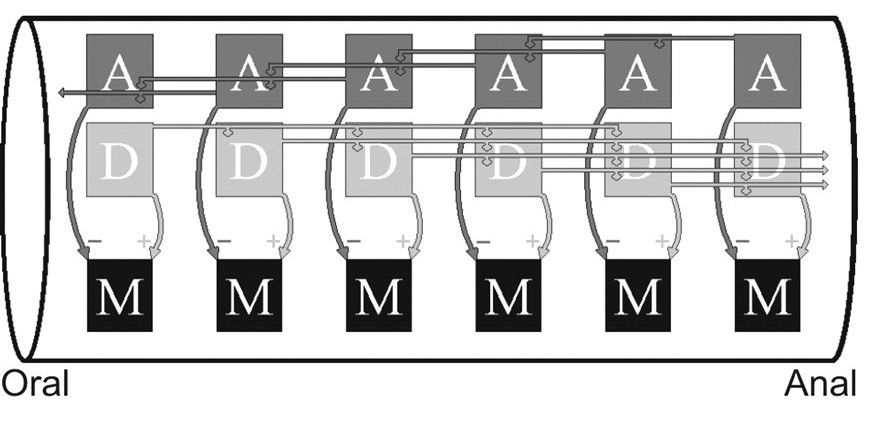}
}
\caption{Schematic `cartoon' of the model of  \cite{KF99}. Each module
(from arbitrarily many) has ascending interneurons and associated motor neurons (A), descending interneurons and associated motor neurons (D), and circular muscle (M). Small gray arrows represent synaptic couplings.}
\label{F:KF99fig}
\end{figure}
\end{example}

\begin{figure}[h!]
\centerline{%
\includegraphics[width=.7\textwidth]{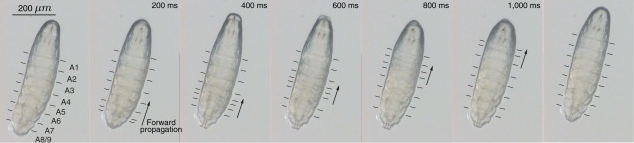}
}
\vspace{.1in}
\centerline{%
\includegraphics[width=.5\textwidth]{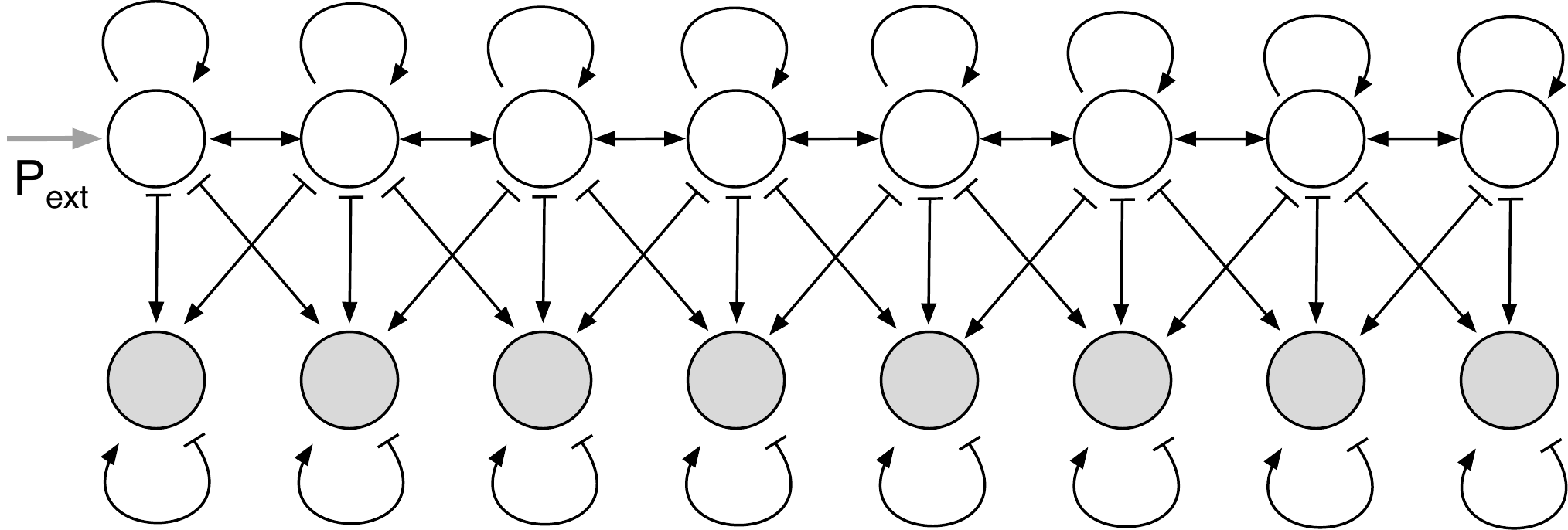}
}
\caption{ {\em Top}: Peristaltic waves of contraction during forward 
 crawling in a first instar {\em Drosophila} larva, courtesy of Gjorgjieva \etal~{\rm \cite{GBEE13}}. {\em Bottom} Model network.}
\label{F:Gjorg1}
\end{figure}

\begin{example}\em
\label{ex:drosophila}

Gjorgjieva \etal~\cite{GBEE13} study neural networks for crawling movement in {\em Dro\-soph\-ila} larvae, which is driven by a peristaltic wave propagating
 from the rear (posterior) to the front (anterior).
Figure~\ref{F:Gjorg1} (top) shows these contractions in the larva, in snapshots taken every 200 ms;
A1--A8/9 indicate the segments;
arrows illustrate simultaneous contraction of neighbouring segments; lines across the larva show the dentical belts, which approximate segment boundaries.
Figure~\ref{F:Gjorg1} (bottom) 
shows the model network studied in that paper. The equations for the dynamics
are of Wilson--Cowan (rate model) type \cite{ET10,WC72}.
The segments are connected 
with nearest-neighbor excitatory connections (triangular arrowheads) 
and inhibitory connections (barred arrowheads). Forward waves are initiated by providing a
time-varying  
external input $P_{\mathrm{ext}}$ into the excitatory population of segment A8.
When $P_{\mathrm{ext}}$ is a short pulse of suitable amplitude a single 
wave is excited; longer pulses excite more waves.

\end{example}

Many networks in the literature have a similar repetitive feedforward structure; 
see for example \cite{MHKDA03,SMG18}. These biological examples are based on
networks of neurons, which control muscle groups, but the general theory applies
more widely. A standard evolutionary pathway is to make multiple copies of
an existing structure, 
and to modify the result through adaptation to different environments.

\subsection{Feedforward Propagation}
\label{SS:FP}
In general, firing signals can propagate naturally along chains of neurons if each 
neuron sends an excitatory signal to the next. However, these signals can lose synchrony
with each other, or phase relations can change, because of random time delays
or other accumulating differences between distinct chains. Similar remarks apply to
other areas of application. A more robust way to propagate dynamic patterns,
in general networks, has its merits.

To set up such a propagation method,
we work in a general context for network dynamics introduced in~\cite{GST05, SGP03},
with slight modifications in \cite{GS23},
which provides a formal framework for analysing networks of coupled
dynamical systems (ODEs). These can be viewed as directed networks
in which nodes and directed edges (`arrows') are labelled
with `types'. Nodes of the same type have the same state space,
and arrows of the same type represent identical types of coupling.
Nodes with isomorphic sets of input arrows obey identical ODEs
when corresponding couplings are identified. We outline this formalism
in Section~\ref{S:NAO}.

The main object of this paper is to use this formalism to
construct, for any small CPG network, a larger network in which the dynamics
of the CPG can be transmitted 
synchronously along chains, trees, or other feedforward cascades of modules. 
The same construction, applied to a CPG with cyclic group symmetry, can lead 
to stable propagation of signals with `phase synchrony' --- identical 
waveforms except for regular phase shifts. We call such behaviour
a {\em phase pattern}. Cyclic group symmetry is intimately involved
in such patterns \cite{GRW12,S20overdet,SP08} and \cite[Chapter 17]{GS23};
see Section~\ref{S:RPPCGS}.
In the context of a chain of successive nodes, 
such states can be viewed as travelling waves.

In this construction the nodes of the modules correspond
to, and have the same types, as the nodes of the CPG --- or, more 
generally, some subset of the CPG. Moreover, any specified
synchrony or phase pattern on the CPG can be extended to the new modules.
Their inputs also correspond to input arrows within the CPG, except that
the tail node for an arrow may be any copy of the
corresponding tail node in the CPG that lies further back along the 
feedforward cascade. This structure implies that any
dynamical state of the CPG (or a subset) can be `lifted' to the entire
cascade by requiring corresponding nodes to be synchronous.
That is, the CPG is a quotient network of the cascade in the sense
of~\cite{GS23,GST05, SGP03}. Conversely, the cascade is a lift of the CPG, so the 
dynamics of the CPG lifts to the feedforward network,
and the modules copy the CPG dynamics.

\subsection{Summary of Paper}

Section~\ref{S:NAO} summarises basic concepts and theorems in the
formalism for network dynamics employed here, with particular emphasis
on balanced colorings, quotient networks, and associated lifts.
We introduce a running example: a 7-node network in which a
directed ring of 3 nodes feed forward into a 4-node chain
as in Figure \ref{F:7nodeFFZ3} of Section \ref{S:AME}. (The numbers 3 and 4 
are chosen for convenience and similar remarks apply for any two positive integers.)

Section~\ref{S:FFL} defines feedforward lifts and establishes their main properties,
especially in the construction of synchronous patterns. We show that the Jacobian 
(derivative) has a
block-triangular form, and use the 7-node example to illustrate this result.

Section \ref{S:BStab} reviews various notions of stability, and the relations between them,
for equilibria and periodic orbits. In particular we discuss Liapunov stability, asymptotic stability,
linear stability, hyperbolicity, and Floquet theory for periodic orbits. We also discuss
a convenient choice of norms for network dynamics.

Section \ref{S:TS} deals with analogous `transverse' stability notions for synchrony-breaking
perturbations of a periodic orbit $\{a(t)\}$ on a CPG network $\GG$ giving rise
to a lifted periodic orbit $\{\tilde a(t)\}$ on a feedforward lift $\widetilde \GG$.
In Theorem~\ref{T:FFStab} we use the block-triangular structure of a feedforward lift
to provide a necessary and sufficient condition for a lifted periodic state 
to be Floquet stable, hence asymptotically stable. This condition
is stated in terms of `transverse Floquet multipliers', which are analogous to the
Floquet multipliers for smaller dynamical systems based on the internal
dynamics of individual nodes. Moreover, we need consider only
the nodes of the CPG. Theorem \ref{T:FFLSstab} provides a similar result
for transverse Liapunov stability.

Section~\ref{S:RTJ} defines the similar but different
condition of `transverse stability' of a synchrony subspace, and relates this to
the diagonal entries of the Jacobian for the CPG. Theorem \ref{T:tstab} shows that
transverse stability implies Floquet stability
for equilibria, and for periodic orbits when node spaces 
are $1$-dimensional. A famous example of \cite{MY60} shows that
this can be false for the periodic case when node spaces have dimension
greater than $1$. We briefly discuss additional conditions that
avoid this problem, together with a related issue: synchronisation of chaotic states. 
This involves a more general concept:
`transverse stability on average' Again this is more satisfactory when
node spaces are $1$-dimensional, and even then,
some aspects are conjectural.

Section~\ref{S:PTW} generalises the transverse Floquet stability condition to
phase-synchronous travelling waves, using a feedforward lift
whose CPG has cyclic group symmetry $\Z_k$.
General results in network dynamics imply that such a network can support
states with phase synchrony, in which the phase shifts are integers multiples of
$T/k$ where $T$ is the overall period \cite[Chapter 3]{GS02}. In a feedforward lift,
these states can be viewed as travelling waves.

Finally, we summarise the main conclusions in Section~\ref{S:C}. 

\section{Networks and Admissible ODEs}
\label{S:NAO}

We briefly review the formalism for network dynamics of~\cite{GST05, SGP03},
taking into account minor improvements introduced in the monograph \cite{GS23}. 

A {\em network} is a directed graph whose nodes and directed edges (`arrows')
are classified into types. A {\em node space} --- usually a finite-dimensional
real vector space --- is assigned to each node, defining a {\em node variable},
which may be multidimensional. The network then
encodes a class of {\em admissible ODEs}, coupled in a manner that respects the
network topology and the node- and edge-types. We give an example in Section~\ref{S:AME}
and a precise definition in Section~\ref{S:AMGC}.

The {\em nodes} (previously called `cells', a term we avoid because
of potential confusion with biological cells) form a (usually finite) set  $\CC = \{1, 2, \ldots, n\}$,
connected by a set $\EE$ of {\em arrows}. 
Each arrow $e$ has a {\em head node} $\hd(e)$ and
a {\em tail node} $\tl(e)$. Nodes are classified into {\em node-types}, and in the associated
admissible ODEs, nodes of the same type have the same internal dynamic.
They also have the same state space, but this property is best treated separately
using the notion of {\em state type}, \cite[Section 9.3]{GS23}. 
Arrows are also classified into {\em arrow-types}, and arrows of the same type
determine the same coupling structure. 

\begin{remark}
In contrast to the conventions in some areas of application
where there are standard model ODEs,
the network diagram does not encode a specific ODE (subject perhaps
to choices of parameters such as reaction rates), and
individual nodes and arrows
do not correspond to specific {\em terms} in a model ODE.
Instead, the network diagram encodes the class of {\em all} ODEs whose 
couplings model the network architecture. This convention is chosen for
mathematical reasons, notably generality \cite[Section 8.10]{GS23}.
\end{remark}

\subsection{Network Diagrams}

A network can be represented graphically by its {\em diagram}, which
is an elaboration of a directed graph. In graph-theoretic
terms it is a {\em coloured digraph}, with colours of nodes and edges 
to represent their node-types, but we use colours in a different manner
so we avoid this terminology. Instead, 
nodes are drawn as dots, circles, squares, hexagons, and so on, with a different
symbol for each type; arrows are similarly decorated to
distinguish arrow-types by using dotted or wavy lines, different shapes of arrowhead,
and so on. Each arrow $e$ runs from the {\em tail node} $\tl(e)$ to the {\em head node} $\hd(e)$. 

An arrow can have the same head and tail, forming a
{\em self-loop}. (A biological term is `autoregulation'.) Two distinct arrows
can have the same head and the same tail, giving
{\em multiple arrows} between the two nodes. This convention
is motivated by some applications and by a basic theoretical
construction, the `quotient network', related to synchrony; see Section~\ref{S:QNL}.

A network $\GG'$ is a {\em subnetwork} of $\GG$ if the nodes of $\GG'$
are a subset $\CC' \subseteq \CC$ and the arrows of $\GG'$ are precisely those
of $\GG$ whose head and tail both lie in $\CC'$.

\subsection{State Spaces}

In order to set up an ODE, we must choose its variables, and the
functions that determine their derivatives. In dynamical systems theory 
the variables determine points in the {\em state space} or {\em phase space}
of the system, which is usually a manifold or more generally a metric space.
Because the term `phase' has other meanings in dynamics, we prefer
the former term. For each node $c \in \CC$, choose
a {\em node (state) space} $P_c$. 
In general, this can be a smooth manifold, and the basic theory of
admissible ODEs and quotient networks remains valid in this context \cite{AF10a,AF10b}.
For simplicity we follow \cite{GS23,GST05,SGP03} and
assume that $P_c = \R^{n_c}$ is a real vector space. (This assumption 
is sufficient for local bifurcation analysis, even if node spaces are manifolds.)
Systems of {\em phase oscillators}~\cite{K88,KE88,KE90,K84}, another standard choice,
correspond to $P_c = \Sone$, the circle.

Phenomena such as synchrony require comparison between
distinct node variables, and this makes sense only when the corresponding
state spaces are equal. State types encode this information:
if nodes $c,d$ are state-equivalent then we require
$P_c = P_d$. The {\em total state space} of the network is the direct sum
\[
P = \bigoplus_{c\in\CC} P_c
\]
and a state is represented by a vector 
\[
x = (x_c)_{c\in\CC}
\]
The entries $x_c$ are themselves vectors when $n_c>1$.

\subsection{Input Sets}

The dynamics of a node depends on the dynamics of its inputs.
We therefore define the {\em input set} of node $c$ to be the set
$I(c)$ of all arrows $e$ such that $\hd(e) = c$. Arrows are used here because
networks can have self-loops and multiple arrows, so specifying the
head and tail does not single out a unique input arrow.

An {\em input isomorphism} $\beta:I(c) \rightarrow I(d)$ is a one-to-one correspondence
between their input sets that preserves arrow-type. That is, $e$ has the same
arrow-type as $\beta(e)$ for all $\beta$ and all $e \in I(c)$
Nodes $c, d$ are {\em input isomorphic} 
if there exists an input isomorphism $\beta:I(c) \rightarrow I(d)$.  Equivalently,
$c$ and $d$ have the same node-type and the same number of input arrows of
each arrow-type.

\subsection{Admissible Maps: Example}
\label{S:AME}

To each network $\GG$ and choice of node spaces $P_c = \R^{n_c}$, 
we associate the class of all ODEs that are
compatible with the network architecture. Such ODEs are called
{\em network ODEs} (previously {\em coupled cell systems}). 
They are determined by the space of {\em admissible vector fields}. 
When all $P_c$ are real vector spaces we refer to these as {\em admissible maps}.
For simplicity we work throughout in the $C^\infty$ category, but most results hold for
$C^r$ with $r \geq 1$.

\begin{example}\em
\label{ex:Z3chain1}

We introduce an example which is revisited several times for different
purposes. Figure~\ref{F:7nodeFFZ3} is a 7-node network, forming
a feedforward chain with a single feedback connection from
node 3 to node 1. (Later, nodes $\{1,2,3\}$ and connecting
arrows are interpreted as a CPG with $\Z_3$ symmetry, 
and the rest of the network is a feedforward lift.)
There is one state-type (all nodes have the same node space), 
one node-type (all nodes have the same type of internal dynamic),
 and one-arrow type (all couplings are identical in form but
 relate to different pairs of nodes). The `colours' of the nodes 
(white, grey, black) are explained in Section~\ref{S:SBC}
and can be ignored here. 

The numbers $3$ and $7$ are
for purposes of illustration, and have no special significance apart
from convenience. Similar examples can be constructed for
any positive integers $p < q$.

\begin{figure}[h!]
\centerline{%
\includegraphics[width=0.6\textwidth]{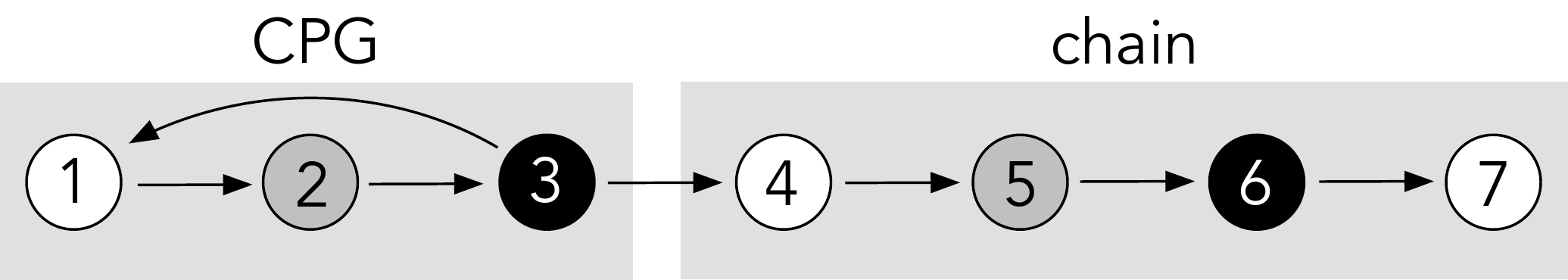}
}
\caption{A 7-node feedforward chain with
one node-type and one-arrow type. Colours show a synchrony pattern.}
\label{F:7nodeFFZ3}
\end{figure}

Admissible ODEs for this network have the following form:
\begin{equation}
\label{E:7nodeFFZ3}
\begin{array}{rcl}
\dot{x}_1 &=& f(x_1,x_3) \\
\dot{x}_2 &=& f(x_2,x_1) \\
\dot{x}_3 &=& f(x_3,x_2) \\
\dot{x}_4 &=& f(x_4,x_3) \\
\dot{x}_5 &=& f(x_5,x_4) \\
\dot{x}_6 &=& f(x_6,x_5) \\
\dot{x}_7 &=& f(x_7,x_6) 
\end{array}
\end{equation}
The same function $f$ is used for all components because all 
nodes have the same node-type and all arrows
have the same arrow-type. The $x_c$ for $1 \leq c \leq 7$
all belong to the same node space $\R^k$ because all nodes
have the same state type. (Indeed, using the same $f$ throughout
requires all node spaces to be the same.) The component
for node $c$ is $f(x_c,x_{i(c)})$ where $i(c)$ is the tail of
the (here unique) input arrow to $c$. In this manner, the admissible
ODEs are precisely those that respect the network structure, 
including preserving node- and arrow-types.
\end{example}

\subsection{Admissible Maps: General Case}
\label{S:AMGC}

We now describe, informally, a procedure for writing down admissible maps
for general networks. Formal definitions are given in 
\cite[Section 9.4]{GS23} and \cite[Section 3]{GST05}. 

For each node $c$ choose {\em node coordinates} $x_c$ on $P_c$.
Nodes of the same state-type have the same coordinate system.
In general, $x_c$ may be multidimensional ($n_c > 1$). Let $P = \oplus_c P_c$
be the total state space.
A map $f = (f_1, \ldots, f_n)$ from $P$ to itself has components
\[
f_c: P \rightarrow P_c \qquad 1 \leq c \leq n
\]
For admissibility we impose extra conditions on the $f_c$ that reflect network
architecture, as follows:

\begin{definition}\em
\label{D:admiss}
Let $\GG$ be a network. 
A map $f: P \rightarrow P$ is $\GG$-{\em admissible} if:
\begin{itemize}
\item[\rm (1)] {\em Domain Condition}:  For  every node $c$, the component $f_c$ 
depends only on the node variable $x_c$ and the input variables $x_{\tl(e)}$ where
$e \in I(c)$.
\item[\rm (2)] {\em Symmetry Condition}:  If $c$ is a node, $f_c$ is
invariant under all permutations of tail node coordinates for equivalent input arrows.
\item[\rm (3)] {\em Pullback Condition}:  If nodes $c \neq d$ are input isomorphic,
the components $f_c, f_d$ are identical as functions. The variables to
which they are applied correspond under some (hence any, by condition (2)) 
input isomorphism.
\end{itemize}
\end{definition}

Formally, conditions (2) and (3) are combined into
a single {\em pullback condition} applying to any pair $c, d$ of nodes, 
equal or different \cite[Remarks 9.20]{GST05}.

Each admissible map $f$ determines an {\em admissible ODE}
\begin{equation}
\label{E:admissODE}
\dot x = f(x)
\end{equation}
where the dot indicates the time-derivative.
If $f$ also depends on a (possibly multidimensional)
parameter $\lambda$, and is admissible as a function
of $x$ for any fixed $\lambda$, we have an {\em admissible family} of maps $f(x,\lambda)$
and ODEs $\dot x =f(x,\lambda)$. Such families arise in bifurcation theory.

\subsection{Synchrony and Balanced Colourings}
\label{S:SBC}

Nodes $c,d$ are {\em synchronous} on a solution $x(t)$ of an admissible ODE if
\[
x_c(t) \equiv x_d(t) \quad \forall t \in \R
\]
This equation makes sense only when $P_c = P_d$; that is, $c$ and $d$ have the same state-type.
Patterns of synchrony that arise naturally and robustly for {\em any} admissible ODE
for a given network are characterised by a property known as balance,
which we now define.

\begin{definition}\em
\label{D:balance}
(a) A {\em colouring} of a network $\GG$ is a map $\kappa: \CC \rightarrow \KK$,
where $\KK$ is a finite set of {\em colours}.

(b) Nodes $c, d$ {\em have the same colour} if $\kappa(c) = \kappa(d)$.

(c) The colouring $\kappa$ is {\em balanced} if there exists
a {\em colour-preserving} input isomorphism for any two nodes of the same colour.
That is, whenever nodes
$c, d$ have the same colour, there exists an input isomorphism
$\beta : I(c) \rightarrow I(d)$ such that $\tl(e)$ and $\tl(\beta(e))$ have the same
colour for all arrows $e \in I(c)$. In symbols, $\kappa(\tl(e))= \kappa(\tl(\beta(e)))$.
\end{definition}
In particular, this definition requires
nodes of the same colour to be input isomorphic. However, the relation
of input isomorphism need not be balanced.

\begin{definition}\em
\label{D:synchspace}
The {\em synchrony subspace} defined by a colouring $\kappa$
of $\GG$ is the vector subspace
\[
\Delta_\kappa = \{x \in P : \kappa(c) = \kappa(d) \implies x_c = x_d \}
\]
That is, nodes of the same colour are synchronous for $x \in \Delta$.
\end{definition}

\begin{example}\em
\label{ex:Z3chain2}
Continuing Example~\ref{ex:Z3chain1},  we again
Consider the 7-node chain of Figure~\ref{F:7nodeFFZ3}. 
The colouring $\kappa$ illustrated in Figure~\ref{F:7nodeFFZ3} has 
three colours $\KK= \{\mbox{B,G,W}\}$, using the initials of the colours black, grey, and white. 
We have
\[
\kappa(1) = \kappa(4) = \kappa(7) = \mathrm{W} \qquad
\kappa(2) = \kappa(5) = \mathrm{G} \qquad
\kappa(3) = \kappa(6) = \mathrm{B} 
\]
All nodes have the same node-type and a single input arrow, and
all arrows have the same arrow-type,
so the nodes are input isomorphic. The colouring is balanced because:
\begin{equation}
\label{E:Z3chain_balance}
\begin{array}{l}
\mbox{Every B node has a single input from a G node.}\\
\mbox{Every G node has a single input from a W node.}\\
\mbox{Every W node has a single input from a B node.}
\end{array}
\end{equation}

All nodes have the same state-type so $P_1 = \ldots = P_7$.
The synchrony subspace is
\begin{equation}
\label{E:Delta_kappa}
\Delta_\kappa = \{(x,y,z,x,y,z,x): x,y,z \in P_1 \}
\end{equation}
\end{example}

The basic theorem on balanced colourings and flow-invariance is:

\begin{theorem}
\label{T:bal_poly}
A subspace $V \subseteq P$  is invariant under every admissible map 
if and only if $V$ is a synchrony space $\Delta_\kappa$ where $\kappa$ is balanced.
\end{theorem}
\begin{proof}
See \cite[Theorem 10.21]{GS23}.
\end{proof}

Theorem~\ref{T:bal_poly} implies that when $\kappa$ is balanced, initial conditions
that have the synchrony pattern defined by $\kappa$ (that is, lie in
$\Delta_\kappa$) give rise to solutions with the same synchrony pattern.
However, this result does not guarantee that the synchrony pattern
is stable: perturbations that break synchrony could cause
the orbit to deviate from $\Delta_\kappa$. This kind of stability 
depends on the admissible map and the orbit concerned. 

\subsection{Quotient Networks and Lifts}
\label{S:QNL}

Balanced colourings give rise to an important construction
in which synchronous nodes are identified in {\em clusters} (or
{\em synchrony classes} or {\em colour classes}).

\begin{definition}\em
\label{D:quot}
Let $\kappa$ be a balanced colouring on a network $\GG$ with colour set $\KK$.
The {\em quotient network} $\GG_\kappa$ has $\KK$ as its
set of nodes (that is, there is one node per colour). 

The node type of node $i \in \KK$ is that of any node $c \in \CC$ such that $\kappa(c) = i$. 

The arrows in $I(i)$ in $\GG_\kappa$ are obtained from
the input set $I(c)$ of any node $c$ with colour $i$ by copying
each arrow $e$ to create an arrow with head $\kappa(\hd(e))$ and
tail $\kappa(\tl(e))$, of the same type as $e$.

The set of arrows of $\GG_\kappa$ is the union of the $I(i)$
as $i$ runs through $\KK$.
\end{definition}

\begin{figure}[h!]
\centerline{%
\includegraphics[width = .15\textwidth]{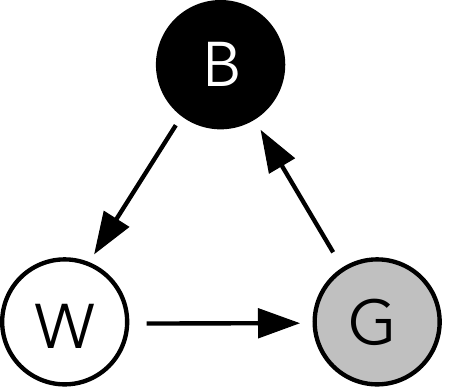}
}
\caption{Quotient network for the balanced colouring of Figure~\ref{F:7nodeFFZ3}
is a unidirectional ring with $\Z_3$ symmetry, permuting nodes and arrows cyclically.}
\label{F:Z3quotring}
\end{figure}

\begin{example}\em
\label{ex:Z3chain3}
The quotient network for the balanced colouring $\kappa$ of 
Figure~\ref{F:7nodeFFZ3} has three nodes $\{\mbox{B,G,W}\}$.
All nodes have the same node type. From \eqref{E:Z3chain_balance}
there is a single arrow-type, with arrows from B to W, W to G, and G to B.
In other words, the quotient network is a $\Z_3$-symmetric unidirectional ring,
Figure~\ref{F:Z3quotring}. In this case
it is isomorphic to the subnetwork $\GG$ with nodes $\{1,2,3\}$
and their connecting arrows. In general, quotient networks need not be subnetworks.
\end{example}

The state space $P_\kappa$ for the quotient network is not the same as $\Delta_\kappa$,
but they can be canonically identified by the map
\begin{equation}
\label{E:nu}
\nu:\Delta_\kappa \to P_\kappa \qquad \nu(x)_{\kappa(c)} = x_c
\end{equation}
which is well-defined.
For example, in \eqref{E:Z3chain_balance}, $\nu(x,y,z,x,y,z,x) = (x,y,z)$.
The projection $\nu$ preserves the synchronous dynamics for any admissible ODE.

\begin{theorem}
\label{T:quot}
Let $\kappa$ be a balanced colouring of $\GG$. Then
\begin{itemize}
\item[\rm (1)]
The restriction of any $\GG$-admissible map to $\Delta_\kappa$
is $\GG_\kappa$-admissible.
\item[\rm (2)]
Every $\GG_\kappa$-admissible map is the restriction to $\Delta_\kappa$
of a $\GG$-admissible map. 
\end{itemize}
\end{theorem}
Another way to say (2) is that every $\GG_\kappa$-admissible map on  $\Delta_\kappa$ {\em lifts}
to a $\GG$-admissible map on $P$.

If $f$ is $\GG$-admissible, the restricted map $f|_{\Delta_\kappa}$ determines the 
dynamics under $f$ of the synchronous clusters determined by the colouring $\kappa$.

Quotient networks can have self-loops and multiple arrows, even if the original
network does not. This feature is required to prove property (2); see \cite[Section 8.10]{GS23}.

\begin{example}\em
\label{ex:Z3chain4}
Again consider the balanced colouring $\kappa$ of 
Figure~\ref{F:7nodeFFZ3}. By Example \ref{ex:Z3chain3} the
quotient network has three nodes $\{\mbox{B,G,W}\}$ forming a $\Z_3$-symmetric ring.
We can write the corresponding coordinates as $x, y, z$ respectively.
Substitute these coordinates, as in \eqref{E:Delta_kappa}, into
the admissible ODE \eqref{E:7nodeFFZ3}:
\begin{equation}
\label{E:7nodeFFZ3_quot} 
\begin{array}{rcl}
\dot{x} &=& f(x,z) \\
\dot{y} &=& f(y,x) \\
\dot{z} &=& f(z,y) \\
\dot{x} &=& f(x,z) \\
\dot{y} &=& f(y,x) \\
\dot{z} &=& f(z,y) \\
\dot{x} &=& f(x,z) 
\end{array}
\end{equation}
This list of equations appears to be overdetermined, because there are
three equations for $\dot x$, and two for each of $\dot y$ and $\dot z$.
However, these equations repeat the same equation three or two times.
(This happens precisely because the colouring is balanced. If it were not,
some equations would disagree with others.) The dynamics therefore reduces to
an ODE with one equation for each coordinate:
\begin{equation}
\label{E:7nodeFFZ3_quot2} 
\begin{array}{rcl}
\dot{x} &=& f(x,z) \\
\dot{y} &=& f(y,x) \\
\dot{z} &=& f(z,y)
\end{array}
\end{equation}
This is the most general admissible ODE on the quotient network,
in accordance with the lifting property.
\end{example}

\section{Feedforward Lifts}
\label{S:FFL}

In this section we define feedforward lifts. We observe that (as is well known)
the Jacobian of any admissible map is block-triangular, with one block for the CPG and separate blocks for each node in the chain. We discuss the construction of
feedforward lifts of a given CPG.
(In alternative terminology \cite{BV02,DL14,MLM20}: the CPG is the 
{\em base} of a {\em graph fibration},
whose {\em fibres} are the synchrony classes of nodes; `feedforward' means that the base 
receives no inputs from the rest of the directed graph.)

We begin by summarising some standard results on feedforward networks; 
see \cite[Chapter 4]{GS23} for proofs. The
usual graph-theoretic term for `feedforward' is {\em acyclic}: no closed directed 
path exists. In dynamics, the term `feedforward' is more common.

By a {\em path} in a network we mean a directed path. In such a path,
each node inputs a signal to the next one, so signals propagate along paths,
but usually they do so without being synchronous or phase-synchronous.

\begin{definition} \rm
(a) Node $q$ is {\em downstream} from node $p$ if there exists a path from $p$ to $q$. 

(b) Node $p$ is {\em upstream} from node $q$ if $q$ is downstream from $p$.
 
(c) Nodes $p$ and $q$ are {\em path equivalent} if node $p$ 
is both upstream and downstream from node $q$. That is, there is a directed path
from $p$ to $q$, and a directed path from $q$ to $p$.

(d) A {\em path component} (or {\em strongly connected component} 
or {\em transitive component}) of a network
is an equivalence class of nodes under path equivalence. It is a subset of nodes
that is maximal subject to the existence of a directed path from any node
in the subset to any other node in the subset. Every network is the disjoint union
of its path components. We use the same term
for the subnetwork obtained by including all arrows between the nodes in 
the equivalence class.

(e) Path component $Q$ is {\em downstream} from path component $P$ 
if there exist a node $p \in P$ and a node $q \in Q$ such that $q$ is downstream from $p$. 
Path component $P$ is {\em upstream} from component $Q$ if $Q$ is downstream from $P$.

(f) The set of all components can be given the structure of
a directed graph whose nodes correspond to the components, with an
arrow from component $C_1$ to component $C_2$
if and only if there exist $c_1 \in C_1, c_2 \in C_2$ with an arrow
from $c_1$ to $c_2$. This graph is the {\em component graph}
or {\em condensation} of the original network, Eppstein~\cite{E16}.
\end{definition}

The following result about the feedforward structure of the component graph is 
well known in the theory of directed graphs \cite{S02}. 
\begin{theorem}
\label{T:PCFF}
The path components are connected in a feedforward
manner; that is, the component graph is acyclic. Moreover, there
is a total order on the nodes that is compatible with the feedforward structure.
\end{theorem}
\begin{proof}
Any easy induction. See for example \cite[Theorem 4.11]{GS23}.
\end{proof}

With a compatible order on the nodes, the Jacobian of any admissible map
is block lower triangular, with the blocks determined by the path components:

\begin{proposition}\label{P:J-lower-block}\rm
The Jacobian\index{Jacobian} matrix of any admissible map 
at any point $x$ is block lower triangular, of the form
\begin{equation}  \label{Jac_block}
J=\Matrix{J_1 &  {\bf 0} &  {\bf 0}    & \cdots & {\bf 0} \\
* & J_2 & {\bf 0} & \cdots &  {\bf 0} \\
* & * & J_3 &  \cdots & {\bf 0} \\
\vdots & \vdots  &\vdots & \ddots & \vdots \\
* & * & * & \cdots  *  & J_m}
\end{equation}
where $J_j$ is the Jacobian matrix of $f$ on the $j$th path component
and each ${\bf 0}$ is a block of zeros of the appropriate size. 
\end{proposition}
\begin{proof}
The value of $f_c(x_c, x_{i_1}, \ldots, x_{i_l})$, where the $i_j$ are
the tails of input arrows to node $c$,
is independent of all $x_d$ for $d>c$.
\end{proof}

The triangular form of~\eqref{Jac_block} implies that the eigenvalues of the Jacobian,
including multiplicity, are determined by the diagonal blocks $J_j$. 
The same is true for the Jordan normal form.

\begin{definition}\em
\label{D:FFL}
Let $\GG$ be a network with a set of nodes $\CC$
and a balanced colouring $\kappa$.
A {\em feedforward lift} of $\GG$ is a network $\widetilde{\GG}$ with
nodes $\widetilde{\CC}\supseteq \CC$ and a balanced colouring $\tilde\kappa$
such that:

(a)
$\GG$ is a subnetwork of $\widetilde{\GG}$.

(b)
Every node $d \in \widetilde{\CC}\setminus \CC$ is 
downstream from some node $c \in \CC$.

(c)
The only loops of $\widetilde{\GG}$ are those that lie in $\GG$. 

(d) 
The colouring 
$\widetilde{\kappa}$ on $\widetilde{\GG}$ has the same set of colours as $\kappa$, 
and restricts to $\kappa$ on $\GG$.
\end{definition}

\begin{example}\em
\label{ex:Z3chain5}
Yet again, consider the 7-node chain of Figure~\ref{F:7nodeFFZ3}, 
in which a CPG with $\Z_3$ symmetry
feeds forward into four additional nodes. Let $\GG$ be the subnetwork
whose nodes are $\{1,2,3\}$ together with the arrows that connect them.
Let $\kappa$ assign different colours to each of these nodes.
Let $\widetilde{\GG}$ be the full 7-node network. Then $\widetilde{\GG}$ 
is a feedforward lift
of $\GG$ for the colouring $\tilde\kappa$ illustrated. It is easy to check
properties (a--d).
\end{example}

There is always at least one balanced colouring on a network $\GG$,
namely the {\em trivial colouring} in which all nodes have distinct colours.

It is easy to see that if $\widetilde{\GG}$ is a feedforward lift of $\GG$ 
and the colouring $\kappa$ is trivial, then
the quotient network $\GG_\kappa$ is isomorphic to $\GG$.

\begin{proposition}
{\rm (a)} Every admissible map $f$ for $\GG$ lifts to an admissible
map $\tilde f$ for $\widetilde{\GG}$.

{\rm (b)} The map $\tilde f$ leaves the synchrony subspace
for $\widetilde{\kappa}$ invariant.

{\rm (c)} The quotient network $\widetilde{\GG}_{\widetilde{\kappa}}$
is isomorphic to the quotient network $\GG_\kappa$.

{\rm (d)} Let $\lambda$ be any balanced colouring of
$\GG$ that is coarser than $\kappa$ (meaning that $\kappa(c)=\kappa(d) \Rightarrow
\lambda(c) =\lambda(d)$)
and let $\widetilde{\GG}$ be a feedforward lift of $\GG$. 
Then $\lambda$ lifts to a balanced colouring $\widetilde{\lambda}$
of $\widetilde{\GG}$ with the same set of colours, and this colouring
is coarser than $\widetilde{\kappa} $.
\end{proposition}
\begin{proof}
These are general properties of lifts; see \cite[Theorem 10.27, Proposition 10.38]{GS23}.
\end{proof}

We call any such $\tilde f$ a {\em synchronous lift} of $f$ with 
{\em pattern of synchrony} $\tilde\kappa$.

By construction, when the colouring on $\CC$ is trivial,
every node in $i \in \widetilde{\CC}\setminus \CC$ has the
same colour as a unique node $c \in \CC$. We denote this node by $[i]$.

\subsection{Path Components of Feedforward Lifts}

Recall that a `path component' of a directed graph is also referred to as a `transitive' or
`strongly connected' component, and consists of a set of nodes that is maximal with
respect to the property that any two
nodes $c,d$ are joined by a directed path of arrows from $c$ to $d$.
In particular, there is also a (possibly different) path from $d$ to $c$.

\begin{lemma}
\label{L:pcFFlift}
Every path component of $\widetilde{\GG}$ is either a path component of $\GG$ 
or a single node of $\widetilde{\CC} \setminus \CC$.
\end{lemma}
\begin{proof}
Nodes $c, d$ are path-equivalent if and only if they lie on a closed loop.
All such loops lie in $\GG$.
\end{proof}

Denote the partial derivative of a function $F$ with respect to a (multidimensional) 
variable $x_c$ by $\DD_c F$. We have:

\begin{corollary}
\label{C:pcFFlift1} 
Let $f$ be admissible for $\GG$ and let $\tilde{f}$ be a lift of $f$
to $\widetilde{\GG}$. Order nodes in a manner that is compatible with
the partial order on the component graph, with nodes $1, \ldots, m$ in $\CC$
and nodes $m+1, \ldots, n$ in $\widetilde{\CC}\setminus \CC$.
Let $J$ be the Jacobian of $f = (f_1, \ldots, f_m)$ on $P_1 \oplus \cdots \oplus P_m$.
Then at any given point the Jacobian $\tilde{J}$ of $\tilde{f}$ is block lower triangular, of the form
\begin{equation}  \label{Jac_block2}
\tilde{J}=\Matrix{J &  {\bf 0} &  {\bf 0}    & \cdots & {\bf 0} \\
* & \DD_{m+1} f_{m+1} & {\bf 0} & \cdots &  {\bf 0} \\
* & * & \DD_{m+2} f_{m+2} & \cdots &  {\bf 0} \\
\vdots & \vdots  &\vdots & \ddots & \vdots \\
* & * & * & \cdots   & \DD_n f_n
}
\end{equation}
evaluated at that point.
\end{corollary}

\begin{corollary}
\label{C:pcFFlift2} 
{\rm (a)} The eigenvalues of $\tilde{J}$ at any point in $P$ 
are those of $J$ together with those of
the $\DD_c f_c$, for $m+1 \leq c \leq n$.

{\rm (b)} At any point in the synchrony space $\Delta_{\kappa}$, 
and for $m+1 \leq c \leq n$, we have
\begin{equation}
\label{E:Dc=D[c]}
\DD_c f_c = \DD_{[c]}f_{[c]}
\end{equation}

where $[c]$ is the unique node $c \in \GG$ such that
$\kappa(c) = \kappa([c])$.
In particular, the eigenvalues of $\DD_c f_c$ are the same as the eigenvalues of 
$\DD_{[c]}f_{[c]}$, when evaluated at the same point.
\end{corollary}

\begin{proof}
(a) This follows from the block-triangular structure.

(b)
The pullback condition and the synchrony pattern induced by $\kappa$
easily imply that $\DD_c f_c = \DD_{[c]}f_{[c]}$.
\end{proof}

\begin{remark}
\label{r:selfloop}
There is a minor complication concerning self-loops. We assume
that all self-loops of the CPG are lifted to feedforward arrows
in $\widetilde{\GG}\setminus \GG$. Thus the
matrix $\DD_c f_c$ is the Jacobian for the internal dynamic on node $c$,
ignoring all inputs from other nodes and all self-loops at $c$ (if any exist).
Now the eigenvalues of $\tilde{J}$ are those of $J$ together with those for the
internal part of each diagonal block of $J$. 
\end{remark}

\begin{example}\em
\label{ex:Z3chain7}
Any admissible map $F(x)$ for the 7-node chain of Figure~\ref{F:7nodeFFZ3} has the form
\eqref{E:7nodeFFZ3}.
Thus the Jacobian at a general point $u =(u_1, \ldots, u_7) \in \R^{7k}$ has the 
block form

\begin{equation}
\label{E:tildeJ7node}
\widetilde{J}|_u = 
\left[
\begin{array}{ccc|c|c|c|c}
f_1(u_1,u_3) & 0 & f_2 (u_1,u_3) & 0 & 0 & 0 & 0 \\
	f_2(u_2,u_1)& f_1(u_2,u_1) & 0 & 0 & 0 & 0 & 0 \\
	0 & f_2(u_3,u_2)& f_1(u_3,u_2) & 0 & 0 & 0 & 0 \\
	\hline
	0 & 0 & f_2(u_4,u_3)& f_1(u_4,u_3) & 0 & 0 & 0 \\
	\hline
	0 & 0 & 0 & f_2(u_5,u_4)& f_1(u_5,u_4) & 0 & 0 \\
	\hline
	0 & 0 & 0 & 0 & f_2(u_6,u_5)& f_1(u_6,u_5) & 0 \\
	\hline
	0 & 0 & 0 & 0 & 0 & f_2(u_7,u_6) & f_1(u_7,u_6) 
\end{array}
\right]
\end{equation}

Here we write $f_1, f_2$ for the partial derivatives of $f$ with respect to 
its first and second variables. In our usual notation, $f_i = \mathrm{D}_i f$
for $i = 1, 2$.
The lines indicate the block-triangular structure, with
a $3 \times 3$ block at top left, which we recognise
as the Jacobian of $F$ restricted to $\GG$, and a series of four
blocks $f_1$. These blocks are evaluated at $u$
and need not be equal, but when evaluated at a point in $\Delta_\kappa$ 
they are equal for nodes of the same colour,
by \eqref{E:Dc=D[c]}.

\end{example} 

\subsection{Construction of Feedforward Lifts}

Feedforward lifts of a given network $\GG$ are easy to construct.
Informally, add new nodes one at a time, choosing a colour from those
in the CPG.  Copy the set of input arrows
from the node of this colour in the CPG, wiring each so that its head is the new node
and its tail is any old node with the same colour as the tail of the
corresponding arrow in $\GG$. Repeat.

More formally, let $\GG$ be a network with nodes $\CC = \{1,\ldots, m\}$.
Colour all of its nodes differently, so we can identify the colour set $\KK$ with
$\{1,\ldots, m\}$.
The construction of a feedforward lift with CPG $\GG$ is simple and obvious.
It can be described inductively, one new node at a time. 
Let $\GG_0 = \GG$ and $\CC_0$ = $\CC$. This is a trivial feedforward
lift of $\GG$ with no extra nodes, and starts the induction process.

Assume that $k$ extra nodes have been added, to obtain a
feedforward lift $\GG_k$ of $\GG$ with nodes $\CC_k = \{1,\ldots, m+k\}$.

Add a new node $m+k+1$ to get $\CC_{k+1}$. Assign this
node the same colour as some node $d \in \CC_0$. It remains to define
the input arrows of node $m+k+1$ in a manner that makes the colouring balanced. 
To do so,
copy the input set $I(d)$ via an input isomorphism $\beta$, assigning
all these arrows the new head node $m+k+1$. Now $I(m+k+1) = \beta(I(d))$. 
Rewire the tail node of each arrow $\beta(e) \in I(m+k+1)$
so that its tail $\tl(\beta(e))$ is any node in $\CC_k$ with the same colour as $\tl(e)$.
(This can be $\tl(e)$ itself, but to obtain short-range connections we can use any
 node further along the chain with the required colour.)
Then $\beta$ is a colour-preserving input isomorphism from $I(d)$ to $I(m+k+1)$.
Since all tail nodes of the new arrows lie in $\CC_k$, the resulting network
$\GG_{k+1}$ is a feedforward lift of $\GG$.

For example, in Figure \ref{F:7nodeFFZ3} we have
$\CC_0 = \{1,2,3\}$ and $\GG_0$ is the $\Z_3$-symmetric ring
on those nodes. We want the colouring with colour-classes
$\{1,4,7\}, \{2,5\},\{3,6\}$. To obtain $\GG_1$ we add node $4$, which has the same colour
as node $1$. Node $1$ has a single input arrow with tail node $3$;
copy this arrow so that its head is node $4$, and the tail remains at node $3$
since this is the only node that is earlier than node $4$ in the ordering and has the correct colour.
To get $\GG_2$ add node $5$ and copy the input arrow to node $2$.
This time there are two choices for the tail node: either node $1$ or node $4$.
The figure chooses $4$. To get $\GG_3$ we need an arrow with
head node $6$ and tail node either $2$ or $5$, and similarly for $\GG_4$.
One set of such choices
(with arrows of the shortest possible range) gives Figure~\ref{F:7nodeFFZ3} with the colouring illustrated.

\begin{remarks}\em
(a) As this description makes clear, feedforward lifts are not unique.

(b) A similar construction can be applied when $\GG$ has a nontrivial balanced
colouring $\kappa$. Its description is essentially identical because the only change
is to $\kappa$. Now $\kappa$ lifts to a balanced colouring $\tilde\kappa$ with the same 
set of colours. This is a general property of colourings \cite[Proposition 10.38]{GS23}.
\end{remarks}

\subsection{Notation for Feedforward Lift}
\label{S:NFL}

We use the following notation for a feedforward lift. We choose
a fixed (but arbitrary) balanced colouring $\kappa$ and use this
to construct the feedforward lift $\widetilde\GG$ from a CPG $\GG$.

In general, given a symbol $s$ for an object defined by $\GG$,
we denote its lift by $\tilde s$ and (where appropriate)
the complementary object by $s^*$. 

Thus we denote the CPG network by $\GG$ with nodes $\CC = \{1, \ldots, m\}$.
The feedforward lift is $\widetilde\GG$ with nodes $\widetilde\CC = \{1, \ldots, n\}$,
where $n > m$. We let $\CC^\ast = \{m+1,\ldots,n\}$ be the nodes of the feedforward
chain.

Denote the total state space for $\GG$ by $P$, and that for 
$\widetilde\GG$ by $\widetilde P$. 
For any subset $\QQ \subseteq \widetilde\CC$ let
$P_\QQ = \oplus_{c \in \QQ} P_c$. 

Exceptionally, denote
the node space of node $c$ by $P_c$ for all $c \in \widetilde \CC$,
since this introduces no ambiguity. Similarly, coordinates of $P = P_\CC$
are denoted by $(x_1,\ldots, x_m)$, those on 
$\widetilde P= P_{\widetilde \CC}$ by $(x_1,\ldots, x_n)$, and those on 
$P^* = P_{\CC^*}$ by $(x_{m+1},\ldots, x_n)$. 

If $f:P \to P$ is admissible for $\GG$, its lift is denoted by
$\tilde f: \widetilde P \to \widetilde P$. If $x(t)$ is a solution of
the ODE $\dot x = f(x)$ on $P$, then its lift is $\tilde x(t)$,
and this is a solution of
the ODE $\dot x = \tilde f(x)$ on $\widetilde P$.

Colour all nodes of $\CC$ differently and let
the corresponding synchrony subspace be $\Delta$. With the chosen ordering
of nodes, the natural identification $\nu$ of $\Delta$ with $P_\CC$ in \eqref{E:nu} satisfies
\[
\nu(v) = (v_1, \ldots, v_m) \quad v \in \Delta
\]
Its inverse is $\nu^{-1} (v_1, \ldots, v_m) = V$ where $V_c = v_{[c]}$.
The feedforward structure, combined with the balance condition,
implies that the quotient dynamics on
$\Delta$ identifies with the dynamics of the CPG $\GG$ on $P_\CC$.
That is, the dynamics of $f|_\Delta$ on $\Delta$ is conjugate to that
of $f_\CC$ on $P_\CC$ by the identification $\nu$.

\section{Background on Stability}
\label{S:BStab}

The stability of a state of a dynamical system was defined and analysed
by Liapunov in 1892--93; see \cite{L92}. Several different concepts of
stability are analysed systematically in \cite{BS70} for 
a continuous flow on a metric space. The modern treatment
mainly focuses on flows and diffeomorphisms on smooth (mostly compact) manifolds; it
was initiated by Smale~\cite{S67} and Arnold \cite{A89},
and developed extensively by their students and others. 

We recall some basic concepts related to stability; see for example \cite{HS74}.
We restrict attention to equilibria and periodic orbits.

\subsection{Equilibria}

First we recall four stability notions for equilibria, two of which are equivalent. For further
information see \cite[Chapter 4]{MLS93} and \cite{L64,LL61,L92}.
Let  $x^*$ be an equilibrium point of the ODE 
\begin{equation}
\label{E:usualODE}
\dot x = f(x) \qquad x \in \R^n
\end{equation}
where $f:\R^n \to \R^n$ is smooth (usually we take this to mean $C^\infty$,
but often $C^r$ for $r \geq 1$ suffices).

\subsubsection{Stability Notions for Equilibria}
\label{SS:SNE}

\paragraph{Liapunov Stability}
The notion of Liapunov stability goes back to
Liapunov~\cite{L92} and is the central topic of \cite{LL61}. 
The equilibrium $x^*$ is {\em Liapunov stable} if,
for every $\eps > 0$, there exists $\delta>0$ such that if $\|x(0)-x^*\| < \delta$
then $\|x(t)-x^*\| < \eps$ for all $t > 0$. This notion applies unchanged
to non-autonomous ODEs if an equilibrium exists.

\paragraph{Asymptotic Stability}
The equilibrium $x^*$ is {\em asymptotically stable} if it is Liapunov stable,
and in addition $\delta$ can be chosen so that $\|x(t)-x^*\| \to 0$ as $t \to +\infty$.

\paragraph{Exponential Stability}
The equilibrium $x^*$ is {\em exponentially stable} if there is a neighbourhood $V$
of $x^*$ and constants $K, \alpha > 0$ such that $\|x(t)-x^*\| < K \ee^{-\alpha t}$
for all $x(0) \in V$. 
(For some norm, not necessarily the Euclidean one, we can assume $K=1$.)

\paragraph{Linear Stability}
The equilibrium $x^*$ is {\em linearly stable} if all eigenvalues of the Jacobian 
$\mathrm{D}f$ evaluated at $x^*$ have negative real part.
\vspace{.1in}

Exponential stability implies asymptotic stability, which in turn implies
Liapunov stability. Neither converse is valid in general. Linear stability
is equivalent to exponential stability.

\paragraph{Stable and Unstable Manifolds}
In nonlinear dynamics, emphasis is placed on the concept of hyperbolicity.
An equilibrium $x^*$ is
{\em hyperbolic} if no eigenvalue of $\DD f|_{x^*}$ has zero real part.
The state space $P$ decomposes as a direct sum $P=E^s \oplus E^u$,
where the {\em stable subspace} $E^s$ is the sum of all generalised eigenspaces
for eigenvalues with negative real parts, and the
 {\em unstable subspace} $E^u$ is the sum of all generalised eigenspaces
for eigenvalues with positive real parts.
Near $x^*$ there exists a smooth {\em stable manifold} $W^s$ and an {\em unstable
manifold} $W^u$, tangent repectively to $E^s$ and $E^u$. These manifolds are unique.
If $A$ is not hyperbolic there is also
a {\em centre subspace} $E^c$ with a tangent {\em centre manifold} $W^c$;
in general it is not unique and only $C^k$ smooth.

\subsection{Periodic Orbits}

The theory for periodic orbits is analogous, but there are minor complications.
Consider a periodic orbit $A= \{a(t)\}$ for some initial condition $a(0)$.
Stability concepts for a periodic orbit are generally obtained by
considering a {\em Poincar\'e section} $\Sigma$ transverse to the orbit,
so the orbit intersects $\Sigma $ at a point $a(0)$.
There is a corresponding {\em Poincar\'e map} or {\em first return map}
$\sigma:\widetilde{\Sigma} \to \Sigma$. Here
$\widetilde{\Sigma}$ is a neighbourhood of $a(0)$ such that
$\sigma(\widetilde{\Sigma}) \subseteq \Sigma$. Now $a(0)$
is a fixed point of $\sigma$, and we can think of $\sigma$ as a
{\em discrete} dynamical system of $\Sigma$.

In particular, a periodic orbit $A$ is {\em hyperbolic} if the derivative 
of a Poincar\'e map at the fixed point corresponding to the orbit has no
eigenvalues on the unit circle. (The corresponding Floquet operator
has an eigenvalue 1 corresponding to the direction tangent to the orbit,
but the Poincar\'e map drops the dimension by 1.) 
In general there are stable and unstable subspaces with 
associated smooth tangent manifolds. 
The smoothness properties of the centre manifold are more technical \cite{K67}.

\subsubsection{Floquet Theory}

Classically, the main notion of stability for a periodic orbit $A = \{a(t)\}$
is defined via Floquet theory \cite[Chapter 1 Section 4]{HKW81}. 
Linearise the ODE about the periodic orbit $A$ to obtain a time-dependent ODE 
\begin{equation}
\label{E:Lin}
\dot y = M(t)y
\end{equation}
 where $M(t) = \DD_yf |_{a(t)}$
is $T$-periodic. Floquet's Theorem \cite{F83} states that there is a {\em fundamental
 matrix} $Y(t)$ such that any solution $y(t) = Y(t)v$ for a constant vector $v$.
 Moreover, there exists a $T$-periodic matrix function $P(t)$ and a constant matrix $B$
 such that every fundamental matrix has the form
\begin{equation}
\label{E:compFloq}
Y(t) = P(t)\ee^{Bt}
\end{equation}
See \cite[p.39]{HKW81}.

The eigenvalues $\beta_i$ of $B$ are the {\em Floquet exponents},
and the eigenvalues $\rho_i$ of $\ee^{BT}$ are the {\em Floquet multipliers}.
The matrix $\ee^{BT}$ is uniquely determined by $M(t)$, so the $\rho_i$
are unique. The real parts of the $\beta_i$ are unique, but their imaginary parts
are unique only modulo $2\pi/T$; see \cite[Note 2 p.40]{HKW81}.
The stability condition is that all $\beta_i$ have negative real part
except for a single eigenvalue $0$ given by the orbit itself;
equivalently, all $\rho_i$ lie strictly inside the unit circle except for
a single eigenvalue $1$. The lack of uniqueness does not affect these statements.

In a more modern treatment \cite{GH83} the
matrix $\ee^{BT}$ is essentially the Jacobian of a Poincar\'e map
at the fixed point corresponding to the periodic orbit, reduced by
one dimension to exclude the eigenvalue $1$ along the periodic orbit.

Stability in this sense implies {\em asymptotic stability},
where now we let $A= \{a(t)\}$ for some initial condition $a(0)$;
then there is a neighbourhood $U \supseteq A$
such that if $x(0) \in U$ with orbit $\{x(t)\}$ then 
\begin{equation}
\label{E:asymp_stab}
\lim_{t \to \infty} d(x(t),A) = 0
\end{equation}
where $d(x,A)) = \inf_{a \in A} \|x-a\|$.
Again the convergence is exponential. Moreover,
for each $x(0) \in U$ there exists $\theta \in \R$, depending on $x(0)$, such that
\begin{equation}
\label{E:asymp_theta}
\lim_{t \to \infty} \| x(t) - a(t+\theta)\| = 0
\end{equation}
and the convergence is exponential. See \cite[Theorem (3) p.42]{HKW81}.
The periodic orbit is then said to have {\em asymptotic phase}. 
For later use we state:
\begin{definition}\em
\label{D:isochron}
The submanifold of initial conditions leading to
a given asymptotic phase $\theta$ is called an {\em isochron} \cite{CR95,G75}. The isochrons
fill out a neighbourhood of the stable periodic orbit.
\end{definition}

\begin{remark}
Every Floquet stable periodic orbit is 
hyperbolic. By \cite[Theorem 4.1(f)]{HPS77}, this implies that the orbit persists 
after any sufficiently small $C^1$
perturbation of the ODE (admissible or not), in the sense that
there exists a unique periodic orbit close to the original one. (The theorem is proved
there for a discrete dynamical system, but at the end of the proof it is stated that the result
is also valid for a continuous one.)
This shows that existence and stability of periodic orbits,
deduced from idealised models, persist when the ideal assumptions
are only approximately valid --- provided the approximation is 
close enough. In practice quite large perturbations often preserve existence and stability;
see \cite{SW23} for some numerical experiments on feedforward lifts.
\end{remark}

\subsubsection{Liapunov Stability}

The notion of Liapunov stability transfers to a periodic orbit $\{a(t)\}$ 
via a Poincar\'e map.
More generally, for any orbit $\{a(t)\}$,
define $y(t)=x(t)-a(t)$. Then the non-autonomous ODE 
(called a {\em system of deviations})
\[
\dot y = f(y+a(t))-\dot a(t)
\]
has an equilibrium at $y=0$. The orbit $a(t)$ is defined to be Liapunov stable
if this equilibrium is Liapunov stable.
Since $y(t) = x(t)-a(t)$ we can unravel this definition:

\begin{definition}\em
\label{D:LSorbit}
The orbit $a(t)$ is {\em Liapunov stable}
if, for any $\eps > 0$, there exists $\delta >0$ such that whenever
$\|x(0)-a(0)\| < \delta$ we have $\|x(t)-a(t)\| < \eps$ for all $t \geq 0$.
\end{definition}

\subsection{Remark on Norms}

All norms on $\R^n$ are equivalent, so the stability notions in Section \ref{SS:SNE}
do not depend on the choice of norm. 
In network dynamics a convenient  norm on $\widetilde P$ is
\begin{equation}
\label{E:net_norm}
\|x\| = \|x_1\| + \cdots + \|x_n\|
\end{equation}
where $\|\cdot\|$ is (say) the Euclidean norm on $P_c$, and this is
consistent with \eqref{E:sum_norm}.
In the next section we work with direct sums $\XX\oplus \YY$
of subspaces $\XX,\YY$, and define norms so that
if $\ZZ=(X,Y) \in \XX\oplus \YY$ then
\begin{equation}
\label{E:sum_norm}
\|Z\| = \|X\|+\|Y\|
\end{equation}
This definition is consistent with \eqref{E:net_norm} when
$\XX$ and $\YY$ are sums of node spaces.

\section{Transverse Stability for a Feedforward Lift}
\label{S:TS}

We now come to the central results of this paper. We show that
this type of feedforward synchrony can be very robust
if the node dynamics on the CPG has certain features that are common in models.
Not only is it dynamically stable: it is structurally stable, preserved
when connection strengths, the forms of couplings,
and the dynamical equations for nodes are perturbed slightly.

\subsection{Floquet Exponents for Forced Systems}

We begin with a general result. It is presumably well known,
but we give a proof for completeness.

Let $P=\R^k, Q=\R^l$ and consider a forced ODE (skew-product)
on  $P \oplus Q$:
\begin{eqnarray}
\label{E:dotX} \dot X &=& F(X) \\
\label{E:dotY} \dot Y &=& G(X,Y)
\end{eqnarray}
having a periodic orbit $(X(t),Y(t)) = (a(t),b(t))$ of period $T$.

The linearised ODE around this orbit (that is, the Floquet equation)
is then
\begin{equation}
\label{E:FFlinODE}
\Matrix{\dot U \\ \dot V} = \Matrix{\DD_1F|_{(a(t),b(t))} & 0 \\
		\DD_1G|_{(a(t),b(t))} & \DD_2G|_{(a(t),b(t))}}
\Matrix{ U \\V} \qquad U \in P, V \in Q
\end{equation}
where the notation $\DD_kH|_{c(t)}$ indicates the partial derivative of $H$ with respect to the $k$th variable, considering $c(t)$ as a parameter.

\begin{lemma}
\label{L:FFFloq}
With the above notation, the Floquet multipliers of $(a(t),b(t))$ on $P\oplus Q$
are those of $a(t)$ on $P$, together with those for the time-dependent ODE
\begin{equation}
\label{E:dotVG}
\dot V = \DD_2G|_{(a(t),b(t))}V
\end{equation}

\end{lemma}
\begin{proof}
By the feedforward structure, $a(t)$ is a $T$-periodic orbit of \eqref{E:dotX},
and this is the Floquet equation for $a(t)$ on $P$. The solution of \eqref{E:dotX}
gives the Floquet multipliers for $a(t)$ on $P$.

The subspace $0 \oplus Q$ is invariant under the flow of \eqref{E:FFlinODE},
and when restricted to this space \eqref{E:FFlinODE} becomes \eqref{E:dotVG}.
Since $0 \oplus Q$ is a complement to $P$,
solutions of this equation yield the remaining Floquet multipliers for 
$a(t)$ on $P\oplus Q$.
\end{proof}

\subsection{Transverse Floquet Multipliers and Exponents}

Lemma~\ref{L:FFFloq} implies that
for any feedforward lift of a fixed CPG feeding forward into a
chain with an arbitrary number of nodes,
the computation of Floquet multipliers can be reduced to simple
computations involving only the CPG. Indeed, the rest of the network
need not be chain: the same remark applies to any feedforward lift.

This simplification arises
for two reasons. First, the balanced colouring involved in a feedforward
lift creates multiple eigenvalues of the Floquet operator. Second, 
the feedforward structure of $\widetilde\GG$ induces the block-triangular
structure \eqref{Jac_block2} on the Jacobian, hence on the Floquet equation.
This structure pervades the entire dynamics.

In detail, we first need:

\begin{definition}\em
\label{D:TFE}
(a)
With the above notation, the {\em transverse Floquet equation} for node $c \in \CC^*$ is
\begin{equation}
\label{E:TFE}
\dot y_c = \DD_{[c]}f_{[c]}|_{a(t)} y_c
\end{equation}
where $\DD_{[c]}$ is the partial derivative with respect to $x_{[c]}$.
Observe that this depends only on the ODE for the CPG and the periodic orbit
for those equations.

By Floquet theory, every solution of \eqref{E:TFE} has the form
\[
y_c(t) = P_c(t)\ee^{B_ct}v
\]
for a constant vector $v$. Here $P_c(t)$ is $T$-periodic and $B_c$ is a constant matrix. Then:

(b)
The matrix $M_c = \ee^{B_cT}$ is the {\em transverse Floquet matrix} for node $c$.

(c)
The matrix $B_c$ is the {\em transverse Floquet exponent matrix} for node $c$.

(d)
The periodic orbit $A=\{a(t)\}$ is {\em transversely Floquet stable} at node $c$ if all eigenvalues
of $M_c$ have absolute value $<1$.

Equivalently, all eigenvalues of $B_c$ have negative real part.

(e) The eigenvalues of $M_c$ are the {\em transverse Floquet multipliers} for node $c$.

(f) The eigenvalues of $B_c$ are the {\em transverse Floquet exponents} for node $c$.

(We use the word `the' in (b,c) even though these matrices are not unique, because
 the eigenvalues in (d,e,f) {\em are} unique.)
\end{definition}

\begin{example}\em
\label{ex:7nodeFHNTFE}
We find the transverse Floquet equations for the
network in Figure \ref{F:7nodeFFZ3}, for FitzHugh--Nagumo neurons 
with voltage coupling. The nodes in the CPG
are $\{1,2,3\}$ and the nodes concerned are those in $\CC^* = \{4,5,6,7\}$.

The equations are:
\beqn
\dot V_c &=& V_c(a-V_c)(V_c-1) - W_c + I + \mu V_{c-1} \\
\dot W_c &=& bV_c - \gamma W_c
\eeqn
for $4 \leq c \leq 7$. Here $\mu$ is the coupling strength.
The corresponding diagonal blocks of the Jacobian are
\[
J_c(t) = \Matrix{-3V_c^2+2(a+1)V_c-a & -1 \\ b & -\gamma}
\]
Evaluated at $a(t) = (\alpha(t),\beta(t))$ these become
\[
J_c(t)|_{a(t)} = \Matrix{-3\alpha_{[c]}(t)^2+2(a+1)\alpha_{[c]}(t)-a & -1 \\ b & -\gamma}
\]
which is independent of $\beta$.
Setting $y_c = (v_c,w_c)$
the transverse Floquet equations are
\[
\Matrix{\dot v_c \\ \dot w_c} = J_c(t)|_{a(t)}\Matrix{ v_c \\  w_c}
	= \Matrix{(-3\alpha_{[c]}(t)^2+2(a+1)\alpha_{[c]}(t)c-a)v_c-w_c \\ bv_c-\gamma w_c}
\]
Although $\mu$ does not appear explicitly, it affects the periodic orbit $A$,
and so affects the transverse Floquet equations.
\end{example}

\subsection{Stability Theorem for Feedforward Lift}
\label{S:STFL}

We can now give a sufficient condition for a feedforward lift of an
equilibrium or periodic orbit to be Floquet stable in the full state space; that is, stable
to perturbations that break synchrony as well as those that preserve synchrony.

\begin{theorem}
\label{T:FFStab}
Let $\{\tilde a(t)\}$ be a feedforward lift of the periodic orbit $\{a(t)\}$
on $P_\CC$. Then:

{\rm (a)} The Floquet multipliers for $\{\tilde a(t)\}$ are the
Floquet multipliers for $a(t)$, together with the transverse Floquet multipliers
for all $c \in \CC^*$.

{\rm (b)} The transverse Floquet multipliers for $c \in \CC^*$ are the
same as those for $[c] \in \CC$.

{\rm (c)} $\{\tilde a(t)\}$ is stable on $\widetilde P$ if and only if
$\{a(t)\}$ is stable on $P_\CC$ and, for
 all nodes in $\CC$, all transverse Floquet multipliers 
have absolute value $<1$. 
\end{theorem}

\begin{proof}
Order the nodes so that the CPG has nodes $\CC = \{1, \ldots, m\}$
and the rest of the network has nodes $\CC^*=\{m+1, \ldots, n\}$.
Let $\GG_k$ be the subnetwork with nodes $\{1, \ldots, m+k\}$
together with all arrows linking those nodes. Then $\GG_{k+1}$
is a feedforward lift of $\GG_k$ for $0 \leq k \leq n-m$.

To prove (a) we argue by induction on $k$.
The statement is trivial for $k=0$. 
The step from $k$ to $k+1$ follows from Lemma~\ref{L:FFFloq},
bearing in mind that when $c \in \CC^*$ the domain of $f_c$
is a subspace of $P_1 \oplus \cdots \oplus P_c$, so the 
time-dependent parameter $a(t)$ restricts onto this subspace.

To prove (b), observe that
because the state $\tilde a(t)$ is a lift of $a(t)$, the functions
$f_c(x)$ and $f_{[c]}(x)$ are equal when $x \in \Delta$.
The same holds for their derivatives at points $a(t) \in \Delta$.

Part (c) now follows from (a) and (b).
\end{proof}

This theorem shows that stability of a lifted periodic state depends
only on the ODE for the CPG, and is independent of the number
of nodes in $\CC^*$. Roughly speaking, the full CPG equation
determines the Floquet multipliers for the CPG, and its diagonal
terms determine all the transverse Floquet multipliers, because
these are the same as those for nodes in $\CC$ of the appropriate colour.

The theorem also implies that when two nodes in the feedforward
chain $\CC^*$ are synchronous, their transverse eigenvalues 
are equal. In other words, the Floquet matrix can have
multiple eigenvalues generically,
within the world of network admissible ODEs, even when the network
has no symmetry. This phenomenon is well known for steady states 
(indeed, it happens for the 7-node network); 
feedforward lifts provide a wide range of examples for periodic orbits.

The same goes when two nodes in the feedforward
chain $\CC^*$ are phase-synchronous, since the Floquet multipliers
are invariant under phase shifts; see Theorem \ref{T:TWtranseigen} below.

\subsection{Isochrons}

Recall the definition of an isochron, Definition \ref{D:isochron}.

It is easy to see that when $A$ is stable,
isochrons $I_\theta$ for $A$ in $P$ extend trivially
to isochrons $\tilde I_\theta$ for $\widetilde A$ in $\widetilde P$:

\begin{theorem}
\label{T:isochron_lift}
Let $\pi:\widetilde P \to P$ be projection onto the first $m$ coordinates.
Then $\pi(\tilde I_\theta) = I_\theta$ for any $\theta$.
\end{theorem}
\begin{proof}
By the feedforward structure, any orbit $x(t)$ for $\tilde f$ 
projects to an orbit $\pi(x(t))$ for $f$.
\end{proof}
Using conjugacy by the natural isomorphism $\nu:\Delta \to P_\CC$, 
we obtain a related projection onto
isochrons of $\widetilde A \subseteq \Delta$.

\subsection{Transverse Liapunov Stability}

An analogous result to Theorem~\ref{T:FFStab}  can be proved for Liapunov
stability.
We begin with a more general set-up.
Consider a feedforward ODE of the form:
\begin{eqnarray}
\label{E:dotX2} \dot X &=& F(X) \\
\label{E:dotY2} \dot Y &=& G(Y,X)
\end{eqnarray}
on a state space $\ZZ=\XX \oplus \YY$, with coordinates $Z=(X,Y)$.
Choose norms on $\XX$ and on $\YY$, and define the norm on $\XX\oplus\YY$
by $\|Z\| = \|X\|+\|Y\|$. 

Suppose we have a $T$-periodic orbit $Z(t) = R(t)$. Projecting,
we have
\[
Z(t) = (X(t),Y(t)) \qquad  R(t)= (P(t),Q(t))
\]
on this orbit.

By the feedforward structure, 
\begin{equation}
\label{E:Porbit}
P(t)\ \mbox{ is a}\ T\mbox{-periodic orbit of \eqref{E:dotX}}
\end{equation}

\begin{definition}\em
\label{D:LS}
The orbit $R(t)$ is {\em Liapunov stable} on $\ZZ$ if, for any $\eps >0$,
there exists $\delta > 0$, such that:
\[
\| Z(0)-R(0) \| < \delta \implies \| Z(t)-R(t) \| < \eps \qquad (t > 0)
\]
\end{definition}

\begin{definition}\em
\label{D:TLS}
The orbit $R(t)=(P(t),Q(t))$ is {\em transversely Liapunov stable} on 
$\YY$ if, for any $\eps >0$,
there exists $\delta > 0$, such that
\begin{eqnarray}
\label{E:Xeq} \| X(0)-P(0) \| &<& \delta \quad\mbox{\rm and}\\
\label{E:yeq} \| Y(0)-Q(0) \| &<& \delta \quad\mbox{\rm imply}\\
\| Y(t)-Q(t) \| &<& \eps \qquad (t > 0)
\end{eqnarray}
\end{definition}

We now come to the main theorem of this section:

\begin{theorem}
\label{T:LS+TLS=LS}
$ $

The orbit $R(t)$ is Liapunov stable on $\ZZ$  if and only if
the following two conditions hold:

{\rm (a)} The orbit $P(t)$ is Liapunov stable on $\XX$, and

{\rm (b)} The orbit $R(t)$ is transversely Liapunov stable on $\YY$.
\end{theorem}

\begin{proof}

First, suppose that $R(t)$ is Liapunov stable on $\ZZ$.
We prove that (a) and (b) hold. 

To prove (a), let $\eps >0$.
The orbit $X(t)$ is independent of $Y(0)$ by the feedforward
structure. So we can choose initial conditions such that $Y(0)=Q(0)$
without affecting $X(t)$. 
Then 
\[
\| Z(0)-R(0) \| = \|X(0)-P(0) \| + \|Y(0)-Q(0) \| = \|X(0)-P(0) \|
\]
Therefore
\beqn
 \|X(0)-P(0) \| < \delta &\implies& \| Z(0)-R(0) \| < \delta \\
 &\implies& \| Z(t)-R(t) \| < \eps \qquad (t>0) \\
 &\implies& \| X(t)-P(t) \| +  \| Y(t)-Q(t) \| < \eps \quad (t>0)\ \mbox{by \eqref{E:sum_norm}}\\
 &\implies& \| X(t)-P(t) \|  < \eps \qquad (t>0)
\eeqn
This is (a).

Now we prove (b). 
In this case the $Y$-dynamics depends on $X$, so we cannot choose initial conditions.

Again, let  $\eps > 0$.
Since  $R(t)$ is Liapunov stable on $\ZZ$,
 there exists $\delta_1 > 0$ such that
\[
\|Z(0)-R(0)\| < \delta_1 \implies \|Z(t)-R(t)\| < \eps \qquad (t>0)
\]
By part (a) there exists $\delta_2$ such that
\[
\|X(0)-P(0)\| < \delta_2 \implies \|X(t)-P(t)\| < \eps \qquad (t>0)
\]
Let $\delta = \shf \min(\delta_1,\delta_2)$.

Suppose that 
\beqn
\|X(0)-P(0)\| &<& \delta \\
\|Y(0)-Q(0)\| &<& \delta 
\eeqn
Then 
\beqn
\|X(0)-P(0)\| &<& \delta_2 \\
\|Y(0)-Q(0)\| &<& \delta_1 \\
\|Z(0)-R(0)\| &<& 2\delta \leq \delta_1
\eeqn
Therefore
\[
\|Z(t)-R(t)\| < \eps \qquad (t>0 )
\]
But by  \eqref{E:sum_norm},
\[
\|X(t)-P(t)\|+\|Y(t)-Q(t)\| = \|Z(t)-R(t)\| < \eps
\]
so
\[
\|Y(t)-Q(t)\| \leq \|Z(t)-R(t)\| < \eps \qquad (t > 0)
\]
which proves (b).

For the converse, assume conditions (a) and (b) hold.
We must prove that  $R(t)$ is Liapunov stable on $\ZZ$.

Let $\eps > 0$. 
By (a), there exists $\delta_1 > 0$ such that
\[
\| X(0)-P(0) \| < \delta_1 \implies \| X(t)-P(t) \| < \eps/2 \qquad (t > 0)
\]
 By (b), there exists $\delta_2 > 0$ such that
\beqn
\| X(0)-P(0) \| &<& \delta_2\quad \mbox{\rm and}\\
\| Y(0)-Q(0) \| &<& \delta_2\quad \mbox{\rm imply} \\
 \| Y(t)-Q(t) \| &<& \eps/2 \qquad (t > 0)
\eeqn
Let $\delta = \min(\delta_1,\delta_2)$. Then, by \eqref{E:sum_norm},
\[
\|Z(0)-R(0) \| = \|X(0)-P(0) \|+\|Y(0)-Q(0) \|
\]
so
\beqn
\|X(0)-P(0) \|&\leq& \|Z(0)-R(0) \| \\
\|Y(0)-Q(0) \|&\leq& \|Z(0)-R(0) \| 
\eeqn
Thus if $\|Z(0)=P(0)\| < \delta$ then $\|X(0)-P(0) \|<\delta_1$ and $\|Y(0)-Q(0) \|<\delta_2$.
Therefore
\beqn
\| X(t)-P(t) \| &<&  \eps/2 \qquad (t > 0) \\
\| Y(t)-Q(t) \|& <& \eps/2 \qquad (t > 0)
\eeqn
and
\beqn
\|Z(t)-R(t) \| &=& \|X(t)-P(t) \|+\|Y(t)-Q(t) \| \\
	&<& \eps/2+\eps/2 = \eps
\eeqn
so $R(t)$ is Liapunov stable on $\ZZ$.

\end{proof}

\subsection{Transverse Liapunov Stability for Feedforward Lift}

Now consider a feedforward lift with CPG $\GG$, whose nodes
$\CC$ are $\{1, \ldots, m\}$ and chain $\CC^* = \{m+1, \ldots, n\}$.
Choose the numbering so that all nodes upstream from
$m+k$ lie in $\GG_{k-1}$. (See \cite[Theorem 4.11]{GS23}.)

For $0 \leq k \leq n-m$ define the subnetwork $\GG_k$
to have nodes $\CC_k=\{1, \ldots, m+k\}$ and all arrows of $\tilde\GG$
that connect them. In particular $\GG_0 = \GG$.
All nodes upstream from $\GG_k$ lie in $\GG_{k-1}$ for $k \geq 1$.

Define
\[
\tilde a^k(t) = (\tilde a_1(t), \ldots, \tilde a_{m+k}(t))
\]
We now define transverse Liapunov stability for a feedforward lift:

\begin{definition}\em
\label{D:TLS_FFL}
The lifted periodic orbit $\tilde a(t)$ is {\em transversely Liapunov stable} 
at node $m+k \in \CC^* = \{m+1, \ldots, n\}$ if
$\tilde a^k(t)$ is transversely Liapunov stable for the decomposition
\[
\XX = P_{\GG_{k-1}} \qquad \YY = P_{m+k}
\]
Since $\tilde a_c(t) \equiv a_{[c]}(t)$, we can write this explicitly:
 for all $\eps > 0$ there exists $\delta_{m+k}$ such that
\begin{equation}
\label{E:TLScond}
\begin{array}{l}
\|x_{\CC_{k-1}}(0)-a_{\CC_{k-1}}(0)\| < \delta_{m+k}\quad \mbox{and}\quad
	\|\tilde x_k(0)-a_{[k]}(0)\| < \delta_{m+k}\\
\qquad	\implies \|\tilde x_k(t)-a_{[k]}(t)\| < \eps\quad (t > 0)
\end{array}
\end{equation}

The lifted periodic orbit $\tilde a(t)$ is {\em transversely Liapunov stable}
if it is transversely Liapunov stable at every node in $\CC^*$.

\end{definition}

\begin{theorem}
\label{T:FFLSstab}
The lifted periodic orbit $\tilde a(t)$ is Liapunov stable if and only if
$a(t)$ is Liapunov stable on $P_\GG$ and 
$\tilde a(t)$ is transversely Liapunov stable at every node of $\CC^*$.
\end{theorem}

\begin{proof}
Use induction on $n-m$ and Theorem \ref{T:LS+TLS=LS}.
\end{proof}

\subsection{Relevance to Applications}
\label{S:RA}

We plan to discuss applications of the methods of this paper
in two follow-up papers: one on chains of standard model neurons,
and one on quadruped locomotion. For reasons of space
these applications cannot be included here, but we comment on
one important issue that arises. We give only a
brief sketch; details will appear in those papers.

Let $\kappa$ be the balanced colouring used to construct the feedforward lift.
Unlike transverse Floquet stability, when $\kappa(c)=\kappa(d)$ the
transverse Liapunov stability condition at node $c$ need not be the same as
that at node $d$. This difference arises because \eqref{E:TLScond} 
involves $\tilde x_{\CC_{c-1}}$ at node $c$, but 
$\tilde x_{\CC_{d-1}}$ at node $d$. These need not be equal.
However, in the inductive argument, we assume:
\beqn
&& \tilde x_{c-1}(t)\ \mbox{is close to}\ \tilde a_{c-1}(t) =   a_{[c-1]}(t)\\
&& \tilde x_{d-1}(t)\ \mbox{is close to}\ \tilde a_{d-1}(t) =   a_{[d-1]}(t)
\eeqn
In dynamical systems theory it
 is common to establish Liapunov stability using a Liapunov function \cite{GH83,LL61,L92}.
In the context of feedforward lifts there is a closely analogous
concept of a `transverse Liapunov function'.
In applications, transverse Liapunov stability can sometimes be proved by
constructing a suitable transverse Liapunov function on $P_c$. The
proof that this function has the required properties depends on
estimates on the supremum of $\|a_{[c]}(t)\|$, since this function
is not usually known explicitly. These estimates
remain valid for any trajectory $\tilde x_{c}(t)$ that is close to $a_{[c]}(t)$.
Thus the same transverse Liapunov function can be used for distinct nodes
$c$ and $d$ such that $\kappa(c)=\kappa(d)$. In practice, therefore,
nodes $c$ and $d$ can be dealt with using the same argument.

\subsection{Transverse Asymptotic Stability}

There is an analogous notion of {\em transverse asymptotic stability},
obtained by replacing `Liapunov' by `asymptotic' in the definition. This leads to
a result analogous to Theorem~\ref{T:FFLSstab}. The proof runs along similar 
lines, but is simpler, so we omit it.

\section{Relation to the Transverse Jacobian}
\label{S:RTJ}

The main difficulty when applying Theorem~\ref{T:FFStab}
is the calculation of the Floquet exponents. As remarked in Section \ref{S:intro},
these exponents must be calculated numerically.
Of course, this remains the case for the CPG dynamics alone,
but it is useful to have a general criterion for condition (b) of the theorem
to be valid, even if only heuristically.
We now discuss one approach to this issue.
We use the notation of Section \ref{S:NFL}.

\begin{definition}\em
\label{D:trans_stab}
The synchrony subspace $\Delta_\kappa$ is {\em globally transversely stable}
if for all $c \in \CC$ and all $x \in \Delta_\kappa$, all eigenvalues of 
each diagonal partial derivative $\DD_c f_c|_x$ have negative real part
for all times $t$.

The lifted periodic state $\widetilde A \subseteq \Delta_\kappa$ is {\em transversely stable}
if for all $c \in \CC$, all eigenvalues of 
each diagonal partial derivative $\DD_c f_c|_{a_c(t)}$ have negative real part
for all times $t$.
\end{definition}

If the lift is constructed so that self-loops become feedforward,
the diagonal block $\DD_c f_c|_x$ refers only to the `internal dynamics' of node $c$.
See Remark \ref{r:selfloop}. We can consider
only nodes in $\CC$
because, on $\Delta_\kappa$, all other nodes are synchronous 
with nodes in $\CC$ via the balanced colouring $\kappa$. Therefore
$f_c$ is the same as $f_{[c]}$ on $\Delta_\kappa$.

\subsection{Counterexamples}
\label{S:notFloquet}
Definition \ref{D:trans_stab}  is motivated by the form of \eqref{Jac_block2}. Historically, it was
conjectured for some time that transverse stability for a stable periodic orbit 
implies stability in the usual Floquet sense; see \cite{AZG11,A49}. 
However, despite the terminology, this conjecture is false in general.
The reason is that although a matrix whose eigenvalues all have negative 
real parts is a contraction in some norm 
\cite[Section 9.1 Theorem (a)]{HS74}, the relevant norm can change
along the periodic orbit. In some circumstances this can create a
Floquet multiplier outside the unit circle. We now give two examples
of this phenomenon.

\begin{example}\em
\label{ex:MY}
An explicit instance is the celebrated {\em Markus--Yamabe counterexample}
\cite[Example p.310]{MY60}.
Consider the ODE $\dot x = A(t)x$ on $\R^2$ where
\[
A(t) = \Matrix{-1+\frac{3}{2}\cos^2 t & 1- \frac{3}{2}\sin t \cos t\\
	-1- \frac{3}{2}\sin t \cos t & -1+\frac{3}{2}\sin^2 t
}
\]
For any $t$, the trace of $A(t)$ is $-\shf$ and the determinant is $\shf$, so
the eigenvalues have negative real part. In fact, they are
$\frac{1}{4} (-1\pm \ii \sqrt{7})$,
for any $t$. However, a solution is
\[
x(t) = \ee^{t/2}\Matrix{-\cos t \\ \sin t}
\]
so the zero solution (which is trivially periodic) is unstable.
\end{example}

\begin{example}\em
\label{ex:notFloquet}
A simpler counterexample uses a discontinuous family of maps $M(t)$. 
This family can then be smoothed
without changing the main conclusion. The literature on such `switching'
or `hybrid' systems is extensive: see for example \cite{L03}.

Let $A, B$ be two constant matrices. Define
\[
M(t) = \left\{ \begin{array}{lcl}
B & \mbox{if} & 0 \leq t < 1 \\
A & \mbox{if} & 1 \leq t < 2
\end{array}\right.
\]
and extend periodically to a family of matrices with period $T=2$.

The solution of \eqref{E:Lin} on $[0,2]$ is then
\[
x(t) = \left\{ \begin{array}{lcl}
\ee^{Bt} x(0) & \mbox{if} & 0 \leq t < 1 \\
\ee^{A(t-1)} \ee^{BT} x(0) & \mbox{if} & 1 \leq t < 2
\end{array}\right.
\]
Thus the Floquet operator is $\ee^A\ee^B$.

When the dimension is 1, this equals $\ee^{A+B}$,
but when the dimension is 2 or more and $A$ and $B$ do not commute,
this expression no longer holds. 
The Campbell-Hausdorff formula \cite[V.5 Proposition 1]{J62}
applies instead.

Let
\[
A=\Matrix{-0.5 & 0\\ 2& -0.7} \qquad B = \Matrix{-0.5 & 2 \\0 &-0.7} = A^\mathrm{T}
\]
Both $A$ and $B$ have eigenvalues $-0.5,-0.7 < 0$. Numerically,
\[
\ee^A = \Matrix{0.606 & 0 \\1.099 & 0.496}
\qquad 
\ee^B =  \Matrix{0.606 & 1.099 \\0 & 0.496}
\]
Now
\[
\ee^A\ee^B = \Matrix{0.367 & 0.666 \\ 0.666 & 1.455}
\]
whose eigenvalues are $1.772, 0.051$. The first of these lies
outside the unit circle.

This example can be made smooth by decreasing the
off-diagonal terms of $A$ to zero and then increasing
the off-diagonal term of $B$, over an arbitrarily short interval of time.
the eigenvalues change by an arbitrarily small amount, so the periodic state
remains unstable.
\end{example}

Heuristically, this phenomenon arises because the flow near $A$
travels roughly parallel to $A$, as well as contracting towards $A$
locally in {\em some} norm.
However, the `parallel' flow changes the local norm in which contraction occurs.
The contraction slows down near $A$, while the flow parallel to $A$ 
remains roughly constant, 
and this can prevent overall contraction.

The change in the norm required for the flow to be contracting
is mainly caused by changes in the (generalised)
eigenvectors of the transverse linearised flow. This is why it is not
picked up by the eigenvalues.

\subsection{Equilibria and 1-dimensional Nodes}

This phenomenon does not occur for equilibria. It can also be
avoided in the context of a feedforward lift if the node spaces are
$1$-dimensional. In the statement of this theorem, `stable'
refers to Floquet stability in the periodic case.

\begin{theorem}
\label{T:tstab}
Let $\widetilde\GG$ be a feedforward lift of a network $\GG$.
Let $f$ be an admissible map for $\GG$. Let $A$ be either
an equilibrium, for node spaces of any dimension, or
a periodic orbit for node spaces of dimension $1$. 
Assume that $A$ is stable in $P_\CC$. 
Let $\tilde f$ be the admissible map for $\widetilde\GG$ 
obtained as a lift of $f$, with lifted periodic orbit $\widetilde A$.
If $A$ is transversely stable, then $\widetilde A$ is
stable for $\tilde f$ in $P_{\widetilde\CC}$. 
\end{theorem}

\begin{proof}
The equilibrium case is trivial because the transverse eigenvalues 
are eigenvalues of the Jacobian at the equilibrium point.

The periodic case follows directly from Theorem~\ref{T:FFStab}(b).
It is well known that for a 1-dimensional space the Floquet equation
can be solved analytically; indeed, the (unique) Floquet exponent
is the time-average of the transverse exponent round the periodic orbit.
The argument is so simple we give it here.

Consider a homogeneous linear equation $\dot y = M(t)y$
where $M:\R \to \R$ is $T$- periodic.
The solution for given
$y(0)$ is found by separation of variables, and is
\[
y(t) = \left(\exp \int_0^t M(t)\dd t \right) y(0)
\]
Since $M(t)<0$ for all $t$, we have $\int_0^T M(t)\dd t < 0$.
Thus the transverse Floquet exponent is negative.
\end{proof}

\subsection{Higher-Dimensional Nodes}
\label{S:HDN}

Example~\ref{ex:notFloquet} can be realised in a feedforward lift
without much difficulty. It shows that transverse stability need not imply
Floquet stability (hence asymptotic stability)
when node spaces have dimension greater than $1$.
Now the situation is more delicate. Because transverse stability can often be tackled
analytically, we discuss these issues briefly.

Additional hypotheses can sometimes be used to establish stability.
An extreme case is when all Jacobians $\DD f|_{a(t)}$
have the same eigenspaces. Then we can decompose according to
the eigenspaces and use a uniform estimate on each eigenspace to prove that the flow
is uniformly exponentially contracting in a suitable norm.
More generally, if the Jacobians $\DD f|_{a(t)}$
have approximately the same eigenspaces, in some reasonable sense,
then provided the approximation is sufficiently close, transverse
stability should imply that the lifted state is stable. 

Transverse stability implies that
the trace of the Floquet matrix is negative,
by \cite[Note 3, p.41]{HKW81}. Equivalently, the product of
the Floquet multipliers (CPG and transverse) lies inside the unit circle.
If the CPG is Floquet stable, this implies that the product of
the transverse Floquet multipliers lies inside the unit circle.

\subsection{Synchronisation of Chaotic Signals}
\label{S:SCS}

We digress to discuss analogous issues when
the equilibrium or periodic cycle $A$ is replaced by a chaotic attractor, a setting
widely used in studies of synchronisation of chaotic signals. There is a vast
literature on this topic, in part because of applications to secure communication.
General references include \cite{BPP00,GM04,PRK01}. Theoretical results are presented in \cite{BR97,PC90,PC98,PCJM97}. Applications to communications include
in \cite{C95,PCKHS92,PC95}.

For chaotic states, there are many notions of stability, and the mathematics
is far more technical. A stable chaotic state is an {\em attractor}, but there
are many distinct definitions of this notion \cite{M85}.
Transverse stability for synchronous chaotic dynamics also 
relies on ideas that are to some extent conjectural, such as the existence of
a Sinai-Bowen-Ruelle (SBR or SRB) measure \cite{B71,KH95}. Some of the issues
involved are discussed for discrete dynamics in \cite{ABS94,ABS96}.
Here we resort to a heuristic description because the chaotic case
is a side-issue for this paper --- though an interesting one. 

Suppose that $S$ is an invariant submanifold, $A\subseteq S$, and
$A$ is an attractor for $f|_S$ in $S$, 
for any reasonable definition of `attractor'. 
Let $\mu$ be an invariant measure on $A$. 
Then we might expect $A$ to be an attractor for $f$ provided that
\[
\int_{u\in A} \DD_{c} f_c|_{u}\, d\mu < 0
\]
for all $c \in \CC$. That is, the transverse flow is {\em attracting on average}
near $A$. The hope is that any local expansion is quickly counteracted by a contraction,
and on average the contractions win.

However, the same problem with invariant manifolds of codimension
greater than $1$ occurs. Moreover --- and worse ---
there can be many distinct invariant measures,
including Dirac measures supported on unstable periodic orbits inside $A$.
Stability can also be defined in several ways. If $\mu$ is an SBR
measure, there is a set of points of positive Lebesgue measure
whose averages are determined by the SBR measure. We then
expect almost all (in the sense of Lebesgue measure) initial points near 
$A$ to be attracted to $A$. Some nearby points may be repelled,
but these form a set whose measure tends to zero near the attractor.
Two such behaviours are on-off intermittency \cite{PST93}
and bubbling \cite{ABS94,ABS96}. Also associated with this
set-up is the concept of a riddled basin \cite{AYYK92}. 
These ideas are discussed rigorously in \cite{ABS94,ABS96},
but only for discrete dynamics and an invariant submanifold of
codimension $1$. Even the existence of SBR
measures is itself largely conjectural, proved mainly for Axiom A systems
in the sense of Smale~\cite{S67} and for more recent generalisations \cite{Y02}, although it is
supported by much numerical evidence for other dynamical systems.

\section{Propagation of Travelling Waves}
\label{S:PTW}

We now generalise the setting of Figure~\ref{F:7nodeFFZ3} so that Theorems
\ref{T:FFStab} and \ref{T:tstab}
apply to certain generic classes of discrete rotating wave in a CPG
with cyclic group symmetry, which, as previously remarked, causes the lifted state
to resemble a travelling wave. Moreover, the sufficient condition
can be applied to just one set of orbit representatives in the CPG,
simplifying the calculations involved.

\subsection{Rigid Phase Patterns and Cyclic Group Symmetries}
\label{S:RPPCGS}

Patterns of phase relations in periodic states for network dynamics
are intimately related to cyclic group symmetries, either of the network
or of its quotient by a balanced colouring. This topic
originated in equivariant dynamics \cite{GSS88};
more recent network analogues are discussed comprehensively 
in \cite[Chapter 17]{GS23}. In particular, there are good reasons to suppose that,
subject to some technical conditions, the quotient network
by synchrony must have cyclic group symmetry to
support a discrete rotating wave \cite{GRW12,S20overdet,SP08} in a structurally stable manner.

We summarise some pertinent results.

Suppose that the CPG $\GG$ has a cyclic symmetry group $\Z_n$.
Then the $H/K$ Theorem \cite{BG01,GS02,GS23} implies that there exist admissible
ODEs $\dot x = f(x)$ whose solutions include a {\em discrete
rotating wave} with spatiotemporal symmetry induced from $\Z_n$.
Such states have a `phase shift symmetry' of the form
\[
x_{\alpha(i)}(t) = x_i (t+ kT/n)
\]
where $T$ is the period, $\alpha$ is a generator of $\Z_n$, and $0 \leq k < n$.

When such a state is lifted to $\GG'$ the rotating wave structure
more closely resembles a travelling wave, because the dynamics of successive nodes
along the lifted chain are identical except for a fixed phase shift $kT/n$.

The stability results of Theorems \ref{T:FFStab} and \ref{T:tstab} apply
in particular to such travelling waves. Moreover, the nodes for which we must
check the transverse Floquet exponents and transverse eigenvalues
can be reduced to those in a single set of orbit representatives
for the $\Z_k$-action.

\subsection{Motivating Example}
Consider the 7-node chain $\widetilde{\GG}$ of 
Figure~\ref{F:7nodeFFZ3}. All nodes have the same state-type: let all 
node spaces be $P_c=\R$ so node variables $x_c$
are $1$-dimensional. The network is feedforward except 
for the backward arrow from node $3$ to node $1$. The subnetwork $\GG$ with
nodes $\{1,2,3\}$ and all arrows connecting those nodes 
can be considered as a CPG with $\Z_3$ symmetry, which
feeds forward into the chain $\{4,5,6,7\}$.

Admissible ODEs take the form \eqref{E:7nodeFFZ3}, and
the Jacobian at any point $u =(u_1, \ldots, u_7) \in \R^7$ has the 
block form \eqref{E:tildeJ7node}.

For suitable $f$ the CPG $\GG$ supports a $T$-periodic $\Z_3$ rotating wave of the form
\[
U(t) = (u(t), u(t+T/3), u(t+2T/3))
\]
(or its reversal, which we obtain by replacing $T$ with $-T$).
Lift this periodic state to $\GG'$; as remarked earlier this can be considered
as a travelling wave of the form
\[
(u(t), u(t+T/3), u(t+2T/3),u(t), u(t+T/3), u(t+2T/3),u(t))
\]
The last four diagonal blocks are then
\beqn
B_4(t) = \DD_1f|_{(u(t),u(t+2T/3))} \\
B_5(t) = \DD_1f|_{(u(t+T/3),u(t))} \\
B_6(t) = \DD_1f|_{(u(t+2T/3),u(t+T/3))} \\
B_7(t) = \DD_1f|_{(u(t),u(t+2T/3))} 
\eeqn
As $t$ runs through $[0,T]$, these are all phase-shifted versions of $B_4$.
Indeed, $B_5(t) = B_4(t+T/3), B_6(t) = B_4(t+2T/3), B_7(t)=B_4(t)$.

In particular, if all eigenvalues of $B_4(t)$ have negative real
part on the periodic orbit $\{u(t)\}$, the same holds for $B_5(t),B_6(t)$, and $B_7(t)$.
Theorem~\ref{T:tstab} now implies that the lifted periodic state is stable
provided the rotating wave $\{U(t)\}$ on $\GG$ is stable on $P_1 \times P_2 \times P_3$,
and all eigenvalues of $B_4(t)= \DD_1f|_{(u(t),u(t+2T/3))}$ have negative real part.
Thus the transverse eigenvalues (which here determine 
stability since nodes are $1$-dimensional) 
depend only on the internal dynamic of one node.

This idea generalises to Theorem~\ref{T:TWtranseigen} below.

\subsection{Schematic of Construction}
\label{S:SC}

\begin{figure}[h!]
\centerline{%
(a) \includegraphics[width=0.6\textwidth]{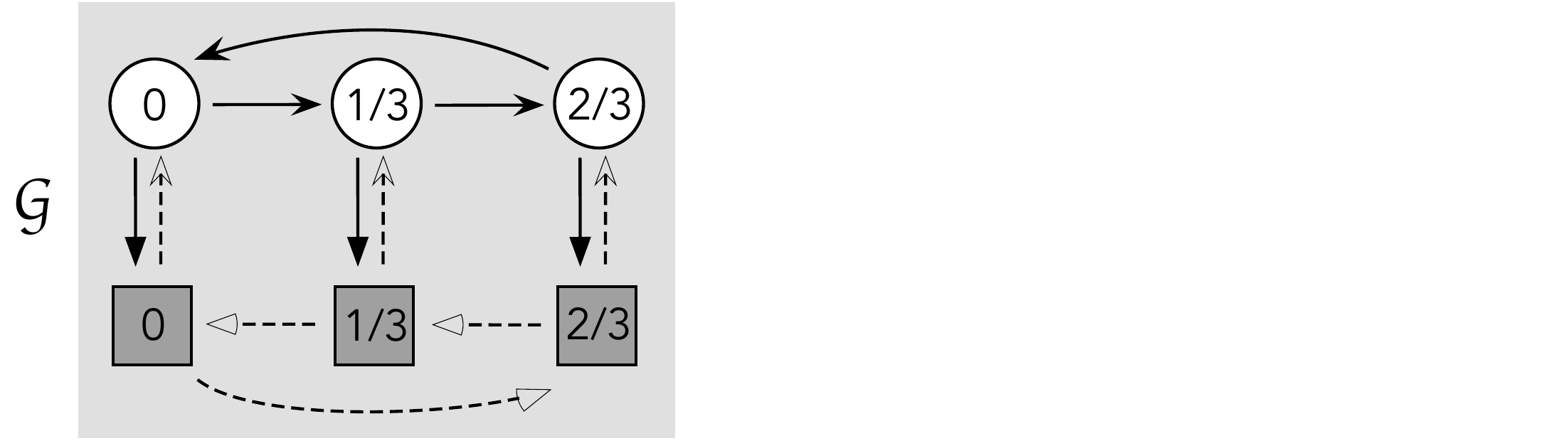}
}
\vspace{.2in}
\centerline{%
(b) \includegraphics[width=0.6\textwidth]{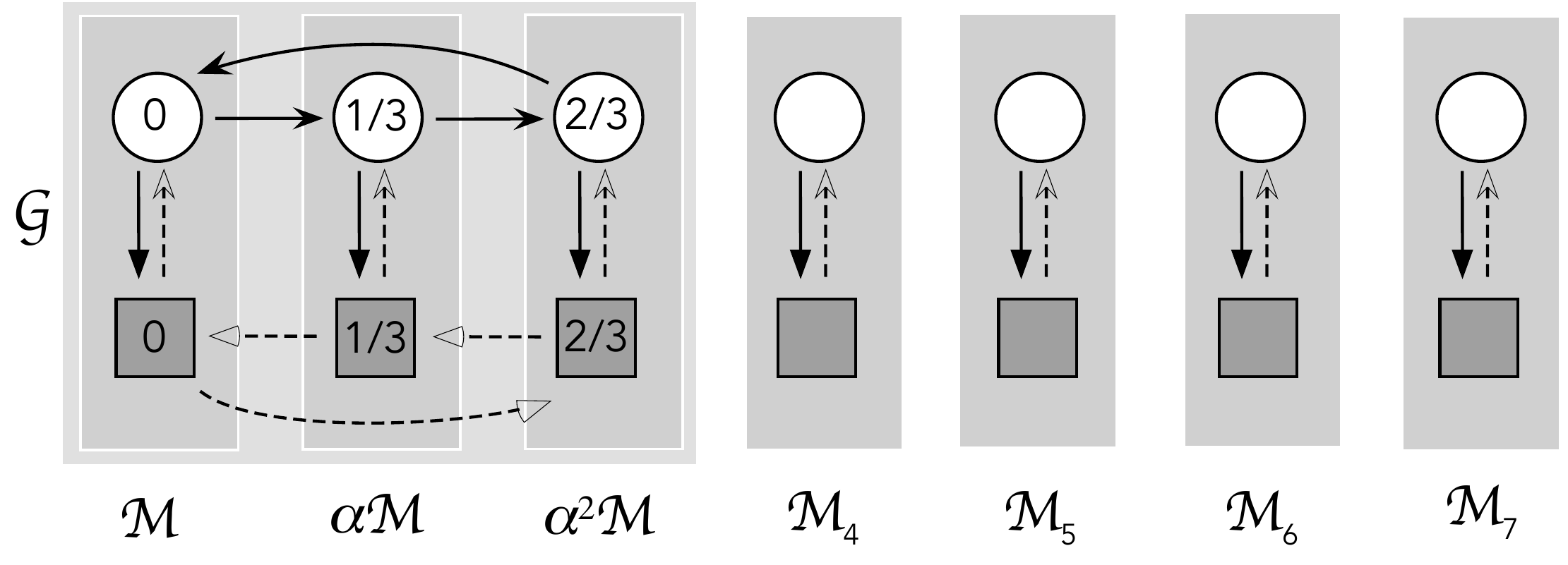}
}
\vspace{.2in}
\centerline{%
(c) \includegraphics[width=0.6\textwidth]{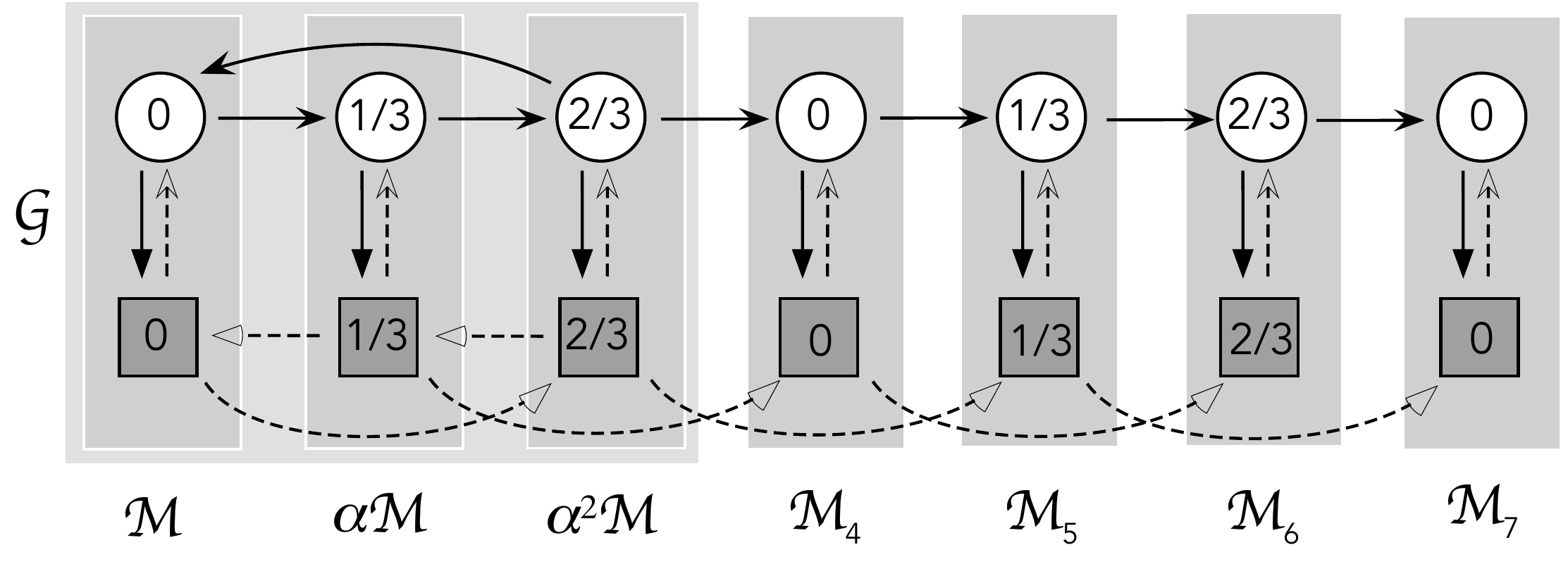}
}
\vspace{.2in}
\centerline{%
(d) \includegraphics[width=0.6\textwidth]{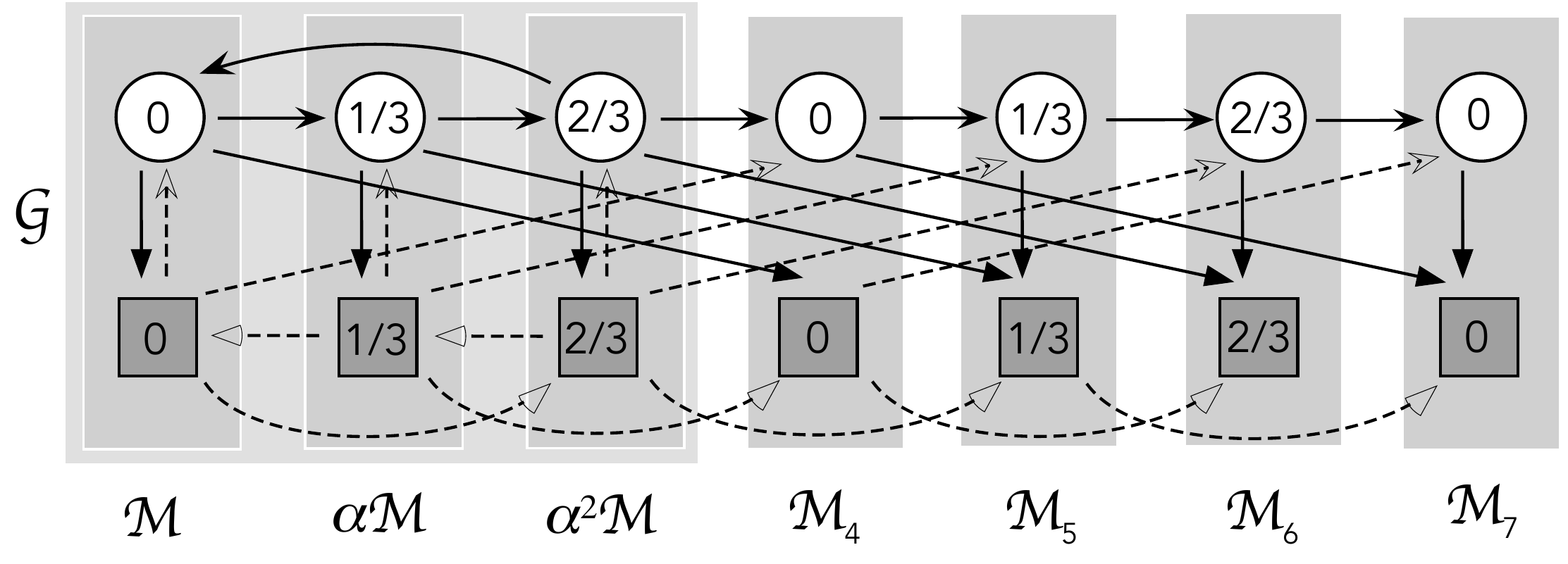}
}
\caption{Four steps in the construction of a feedforward lift for
a CPG with cyclic symmetry supporting a periodic orbit with a rotating wave 
phase pattern. (a) Initial CPG $\GG$ with $\Z_k$ symmetry
(here $k=3$). Fractions $0, 1/3, 2/3$ indicate relative phases.
(b) Module $\MM$, its images under a generator $\alpha$, and
copies  $\MM_4$ -- $\MM_7$. (c) Feedforward inputs to $\MM_4$ -- $\MM_7$. (d) Rewiring
$\MM_4$ -- $\MM_7$ to make all inputs feedforward.}
\label{F:Z3FFchain}
\end{figure}

Figure~\ref{F:Z3FFchain} is a schematic illustration of the four steps involved
in constructing a feedforward lift from a rotating wave state to a travelling wave
along a chain.

(a) Consider a CPG $\GG$ with $\Z_k$ symmetry, generated by
a bijection $\alpha:\CC \to \CC$. (Here $k=3$.) For simplicity, assume
that $\alpha$ is a product of cycles with the same 
length $k$ and every node occurs in one of these cycles.
The set $\CC$ can then be partitioned into $k$ disjoint subsets
that are cycled by $\Z_k$. We refer to any such subset
as a {\em module}, because the same `modular' structure is repeated on each 
$\Z_k$-orbit, modelling the structures in Sections \ref{SS:BM} and \ref{SS:FP}.
Let $f$ be an admissible map such that the ODE $\dot x = f(x)$
has a discrete rotating wave state $u(t)$ satisfying the phase relation

\begin{equation}
\label{E:phaserel}
\alpha u(t) = u(t+T/k)
\end{equation}

(Relative phases marked inside node symbols.)

(b) Consider a module $\MM$; that is, a set of nodes 
comprising one representative from each $\Z_k$-orbit. 

Assign phase $0$ to these nodes, so that the other
 $\Z_k$-orbits correspond to phase shifts $T/k, 2T/k, \ldots,$ $(k-1)T/k$.
 Copy the module (along with any arrows whose heads and tails lie in the module)
 to obtain $\MM_{k+1}, \MM_{k+2}, \ldots, \MM_l$. (Here $l=7$.)

(c) Assign phases $0, T/k, \ldots, (l-k)T/k$ to nodes in  $\MM_{k+1}, \MM_{k+2}, \ldots, \MM_l$.
Assign input arrows to these nodes, preserving the arrow type and the phase
relations in $\GG$. Do so in a manner that makes all new arrows feedforward.

(d) Rewire internal arrows in $\MM_{k+1}, \MM_{k+2}, \ldots, \MM_l$,
preserving the arrow type and the phase
relations in $\GG$. Do so in a manner that makes all rewired arrows feedforward.
(This stage is optional: it simplifies the calculation of Jacobians but may be
less realistic biologically. For example, if modules correspond to segments
of an organism, neuronal connections within segments are likely to be the
same in each segment.) 

\subsection{General Theorem}

We now state a general theorem for such constructions, and 
prove that the transverse eigenvalues depend only on the internal dynamic
of one module.

\begin{theorem}
\label{T:TWtranseigen}
Assume that $\GG$ has nodes $\CC = \{1, \ldots, m\}$ 
with a cyclic automorphism group $\Z_k = \langle \alpha \rangle$,
such that $n=mk$ and $\alpha$ acts like the cycle $(1\, 2\, \ldots\, k)$
on all of its orbits on $\GG$. Let $u(t)$ be a $T$-periodic solution of an admissible ODE
with discrete rotating wave phase pattern~\eqref{E:phaserel}.
Choose a module $\MM$ of orbit representatives.
Let $\widetilde{\GG}$ be obtained
by lifting appropriate copies of translates of this module by $\Z_k$,
as described in Section~{\rm\ref{S:SC}}. Then 

{\rm (a)} The periodic state $u(t)$ on $\GG$ lifts to a $T$-periodic travelling wave
state $\tilde{u}(t)$ for $\widetilde{\GG}$ with phases corresponding to the extra copies
$\MM_{k+1}, \MM_{k+2}, \ldots, \MM_l$ of $\MM$.

{\rm(b)}
The Floquet exponents (evaluated at any point) are 
those on the module $\MM$, together with the transverse
Floquet exponents for $\MM$.

{\rm(c)} If the Floquet exponents  and
the transverse Floquet exponents on $\MM$ have negative real part,
then $\tilde{u}(t)$ is stable.
\end{theorem}

\begin{proof}
(a) This follows because $\GG'$ is a lift of $\GG$.

(b)
Let $\tilde{f}$ be the lift of $f$. The transverse
Floquet exponents must have negative real part for Theorem \ref{T:FFStab}
to apply. The Floquet matrix $\ee^{BT}$ is independent of the initial time
chosen for one period of the flow, hence its eigenvalues are the same after any
phase shift. Therefore phase-synchronous nodes have the same transverse
Floquet exponents.

(c) This follows from Theorem~\ref{T:FFStab}.
\end{proof}

By Theorem~\ref{T:tstab} and the above we immediately deduce:
\begin{corollary}
\label{C:}
{\rm(a)}
The transverse eigenvalues (evaluated at any point) are the same as the eigenvalues of the
Jacobian $J^\MM$ on the module $\MM$, including only the arrows 
whose heads and tails lie in $\MM$ (evaluated at the same point). 

{\rm(b)} If all eigenvalues of $J^\MM$ have negative real part
when evaluated on $\{u(t)\}$, then
$\tilde{u}(t)$ is globally transversely stable.

{\rm(c)} If $u(t)$ is stable on $P^\GG$, nodes are $1$-dimensional,
and all eigenvalues of $J^\MM$ have negative real part
when evaluated on $\{u(t)\}$, then
$\tilde{u}(t)$ is stable on $P^{\widetilde{\GG}}$.
\qed
\end{corollary}

Again we emphasise that transverse stability in this sense applies to the
synchrony subspace $\Delta$, and does not guarantee stability in
the Floquet sense, except when nodes are 1-dimensional or the state
is an equilibrium.

More general results of the same kind can be derived for other
actions of $\Z_k$ on $\GG$, such as those leading to multirhythms
\cite{GNS04,GS23}. In each case the connections in the lift must be
tailored to the phase relations of the periodic state concerned. We
do not state such generalisations but in principle the same ideas apply.

Similar remarks to those in Section~\ref{S:HDN} apply to the
phase-synchronous case.

There is a natural analogue of Theorem~\ref{T:TWtranseigen}(c) for Liapunov
stability, proved in the same manner. We do not state it here.
The remark in Section \ref{S:RA} applies to transverse 
Liapunov stability of phase-shifted periodic states.

\section{Conclusions}
\label{S:C}

Propagation of synchronous or phase-synchronous states along
linear chains is important in biology, medicine, and robotics, among other areas of application.

In this paper we establish a theory of feedforward lifts, which provide
a simple, effective, and robust way to propagate signals
with specific synchrony and phase patterns in a stable manner is to 
use a CPG to generate the underlying patterns and propagate them
along a feedforward chain. Suitably constructed, such a 
feedforward lift preserves the waveform of the signal as it propagates.

Applications of the theory are not included in this paper, for reasons of space.
We plan to discuss them in two follow-up papers: one on chains of standard 
model neurons and another on quadruped locomotion.

An important issue is the stability of the propagating signals. Specifically,
maintaining the synchrony or phase pattern requires stability to synchrony-breaking
perturbations --- transverse stability. We give a necessary and sufficient condition
for stability (in the Floquet sense) that depends only on the internal dynamics of 
the CPG nodes.
This implies that if adding a single copy of the CPG leads to a stable periodic orbit,
the chain can be extended arbitrarily far, and even branch into a tree, with the lifted orbit remaining stable.
Transverse Floquet multipliers for a lifted periodic orbit are generically 
multiple whenever nodes in the lift, but not in the CPG, are synchronous.
Analogous results hold for Liapunov stability.

A simpler condition `transverse stability of the synchrony subspace' 
implies linear stability of equilibria, and Floquet stability of
periodic orbits when nodes are $1$-dimensional. The latter implication can fail for
higher-dimensional nodes, but has some heuristic value.

There is a straightforward generalisation of these results to propagating 
phase patterns, where the CPG is
a symmetric ring of identical modules, and generates a rotating wave
with regularly spaced phase shifts. The lifted periodic orbit can be viewed as
a travelling wave along the chain. Again transverse stability
need be verified only for a single module in the chain.

As a final, more speculative remark: 
Feedforward lifts have a simple modular structure, capable
of generating stable propagating signals with specific phase patterns.
This combination of repetitive modules and potentially useful 
dynamical patterns can evolve
naturally from simpler structures, especially in the context of muscle groups driven
by a network of neurons. This could be one reason why such 
architectures are common in living organisms.

\section*{Acknowledgments}
We thank Peter Ashwin, Marty Golubitsky, and John Guckenheimer for helpful discussions;
two anonymous reviewers for comments that greatly improved the paper; and
and Eddie Nijholt for correcting several typographical errors and noticing that
the proof of Theorem \ref{T:FFLSstab} requires a slightly stronger 
definition of transverse Liapunov stability than the one we originally used.

\bibliographystyle{siamplain}
\bibliography{references}

\begin{thebibliography}{99}


\bibitem{AM78}
R.~Abraham and J.E.~Marsden.  {\em Foundations of Mechanics},
Benjamin/Cummings, New York 1978.

\bibitem{AZG11}
S. Addas-Zanata and B. Gomes.
Horseshoes for a generalized Markus-Yamabe example,
{\em Qualitative Theory Dyn. Sys.} {\bf 10} (2011) 327-332.

\bibitem{AF10a}
N. Agarwal and M.J. Field. Dynamical equivalence of network architecture for coupled dynamical systems I: asymmetric inputs, {\em Nonlinearity} {\bf 23 }(2010) 1245--1268. 

\bibitem{AF10b}
N.  Agarwal and M.J. Field. Dynamical equivalence of network architecture for coupled dynamical systems II: general case, {\em Nonlinearity} {\bf 23 }(2010) 1269--1289.

\bibitem{A49}
M. Aizerman. 
On a problem concerning the stability in the large of a dynamical system,
{\em Uspekhi Mat. Nauk} {\bf 4} (1949) 187--188.

\bibitem{AYYK92}
J.C. Alexander, J.A. Yorke, Z. You, and I. Kan.
Riddled basins, {\em Internat. J. Bif. Chaos} {\bf 2} (1992) 795--813.

\bibitem{A57}
T.M. Apostol.
{\em Mathematical Analysis}, Addison-Wesley, Reading MA 1957.

\bibitem{A89}
V.I. Arnold. {\em Mathematical Methods of Classical Mechanics}, Springer, Berlin 1989.

\bibitem{AP90}
D.K. Arrowsmith and C.M. Place.
{\em An Introduction to Dynamical Systems}, Cambridge University Press, Cambridge 1990.

\bibitem{ABS94}
P. Ashwin, J. Buescu, and I. Stewart.
 Bubbling of attractors and synchronization of chaotic oscillators, 
 {\em Phys. Lett. A} {\bf 193} (1994) 126--139.
 
\bibitem{ABS96}
P. Ashwin, J. Buescu, and I. Stewart.
From attractor to chaotic saddle: a tale of transverse instability, 
{\em Nonlinearity} {\bf 9} (1996) 703--737.

\bibitem{B19}
A. Berkowitz.
Expanding our horizons: central pattern generation in the context of complex activity sequences, 
{\em J. Exp. Biol. } {\bf 222} (2019) 192054;
doi: 10.1242/jeb.192054.

\bibitem{BS70}
N.P. Bhati and G.P. Szeg\"o.
{\em Stability Theory of Dynamical Systems},
Grundlehren {\bf 161}, Springer, Berlin 1970.

\bibitem{BPP00}
S.~Boccaletti, L.M.~Pecora, and A.~Pelaez.
A unifying framework for synchronization of coupled dynamical
systems, {\em Phys.~Rev E} {\bf 63} (2001) 066219.

\bibitem{BV02}
P. Boldi and S. Vigna.
Fibrations of graphs, {\em Discrete Math.} {\bf 243} (2002) 21--66.

\bibitem{B71}
R. Bowen. Periodic points for Axiom A diffeomorphisms,
{\em Trans. Amer. Math. Soc.} {\bf 154} (1971) 377--397.

\bibitem{BBC12}
J.H. Boyle, S. Berri, and N. Cohen.  Gait modulation in {\em C. elegans}:
an integrated neuromechanical model,
 {\em Front. Comput. Neurosci.} {\bf 6} (2012);
doi: 10.3389/fncom.2012.00010.

\bibitem{BR97}
R. Brown and N.F. Rulkov.
Synchronization of chaotic systems: Transverse stability of trajectories in invariant
manifolds, {\em Chaos} {\bf 7} (1997) 395--413; doi: 10.1063/1.166213.

\bibitem{B01}
P.-L. Buono. 
Models of central pattern generators for quadruped locomotion: II.
Secondary gaits, {\em J. Math. Biol.} {\bf 42} (2001) 327--346.

\bibitem{BG01}
P.-L. Buono and M. Golubitsky. 
Models of central pattern generators for quadruped locomotion: I. Primary gaits,
{\em J. Math. Biol.} {\bf 42} (2001) 291--326.

\bibitem{BP04}
P.-L. Buono and A. Palacios. 
A mathematical model of motorneuron dynamics in the heartbeat 
of the leech, {\em Physica D} {\bf 188} (2004) 292--313.

\bibitem{CP83}
R.L. Calabrese and  E. Peterson.
Neural control of heartbeat in the leech {\em Hirudo medicinalis}, 
in: {\em Neural Origin of Rhythmic Movements} (eds. A. Roberts and B. Roberts), 
{\em Symp. Soc. Exp. Biol.} {\bf 37} (1983) 195--221.

\bibitem{CNO95}
R.L. Calabrese, F. Nadim and \O.H. Olsen. 
Heartbeat control in the medicinal leech: A model system for understanding 
the origin, coordination, and modulation of rhythmic motor patterns, 
{\em J. Neurobiol.} {\bf 27} (1995) 390--402.

\bibitem{C95}
T.L. Carroll.
Communicating with use of filtered, synchronized chaotic signals, 
{\em IEEE Trans. Circuits Syst.} {\bf 42} (1995) 105--. 

\bibitem{CBT08}
J.D. Chambers, J.C. Bornstein, and E.A. Thomas. Insights into
mechanisms of intestinal segmentation in guinea
pigs: a combined computational modeling and in
vitro study, {\em Am. J. Physiol. Gastrointest: Liver
Physiol.} {\bf 295} (2008) G534--541.

\bibitem{CTB13}
J.D. Chambers, E.A. Thomas, and
C. Bornstein. Mathematical modelling of enteric neural motor patterns,
{\em Proc. Austral. Physiol. Soc.} (2013)
{\bf 44} 75--84 .

\bibitem{CS93a}
J.J. Collins and I. Stewart.
 Hexapodal gaits and coupled nonlinear oscillator models, 
{\em Biol. Cybern.} {\bf 68} (1993) 287--298.

\bibitem{CS93b}
J.J. Collins and I. Stewart.
Coupled nonlinear oscillators and the symmetries of animal gaits, 
{\em J. Nonlin. Sci.} {\bf 3} (1993) 349--392.

\bibitem{CR95}
S.M. Cox and A.J. Roberts.
 Initial conditions for models of dynamical systems, {\em Physica D} {\bf 85} (1995) 126--141.

 \bibitem{DL14}
L. DeVille and E. Lerman.
Modular dynamical systems on networks, 
{\em J. Eur. Math. Soc.} {\bf 17} (2013);
doi: 10.4171/JEMS/577.

\bibitem{E16}
D. Eppstein.
{\em Design and Analysis of Algorithms}, CS/CSE 161, U. California, Irvine 2016.

\bibitem{ET10}
B. Ermentrout and D. Terman.
{\em The Mathematical Foundations of Neuroscience},
Springer, New York 2010.

\bibitem{F83}
G. Floquet.
Sur les \'equations diff\'erentielles lin\'eaires \`a coefficients p\'eriodiques,
{\em Ann. \'Ecole Norm. Sup. Paris} {\bf 12} (1883) 47--89.

\bibitem{F08}
J.B. Furness. {\em The Enteric Nervous System},
Blackwell, Oxford 2008.

\bibitem{G21}
J.-M. Ginoux.
Slow invariant manifolds of slow-fast dynamical systems,
{\em Internat. J. Bif. Chaos} {\bf 31} (2021) 2150112; 
arXiv:2012.06770.

\bibitem{GBEE13}
J. Gjorgjieva, J. Berni, J.F. Evers, and S.J. Egle.
Neural circuits for peristaltic wave propagation in crawling {\em Drosophila} larvae: analysis and modeling,
{\em Front. Comput. Neurosci.} {\bf 7} (2013); doi: 10.3389/fncom.2013.00024. 

\bibitem{GNS04}
M.~Golubitsky, M.~Nicol, and I.~Stewart. Some curious phenomena in coupled cell networks,
{\em J. Nonlinear Sci.} {\bf 14} (2004) 207--236.

\bibitem{GRW12}
M. Golubitsky, D. Romano, and Y. Wang. Network periodic solutions:
patterns of phase-shift synchrony, {\em Nonlinearity} {\bf 25} (2012) 1045--1074. 
  
\bibitem{GS02}
M.~Golubitsky and I.~Stewart.  {\em The Symmetry Perspective},
{\em Progress in Mathematics} {\bf 200}, Birkh\"auser, Basel 2002.
  
\bibitem{GS23}
M. Golubitsky and I. Stewart. 
{\em Dynamics and Bifurcation in Networks}, SIAM, Philadelphia, to appear 2023.

\bibitem{GSBC98}
M. Golubitsky, I. Stewart, P.-L. Buono, and J.J. Collins.
A modular network for legged locomotion, {\em Physica} D {\bf 115} (1998) 56--72.

\bibitem{GSCB99}
M. Golubitsky, I. Stewart, J.J. Collins, and P.-L. Buono.
Symmetry in locomotor central pattern generators and animal gaits, 
{\em Nature} {\bf 401} (1999) 693--695.

\bibitem{GSS88}
  M.~Golubitsky, I.~Stewart, and D.G.~Schaeffer. {\em Singularities
  and Groups in Bifurcation Theory II}, Applied Mathematics Series,
  {\bf 69}, Springer, New York 1988.

\bibitem{GST05}
M. Golubitsky, I. Stewart, and A. T\"or\"ok.
Patterns of synchrony in coupled cell networks with multiple arrows,
{\em SIAM J. Appl. Dynam. Sys.} {\bf 4} (2005) 78--100.

\bibitem{GM04}
J.M. Gonz\'alez-Miranda. 
{\em Synchronization and Control of Chaos. An introduction for scientists and engineers}, Imperial College Press, London 2004.

\bibitem{Gr03}
H. Gregersen. {\em Biomechanics of the Gastrointestinal Tract}, 
Springer, London 2003.

\bibitem{G75}
J. Guckenheimer.
 Isochrons and phaseless sets, {\em J. Math. Biol.} {\bf 1} (1975) 259--273.

\bibitem{GH83}
J. Guckenheimer and P. Holmes.
{\em Nonlinear Oscillations, Dynamical Systems, and Bifurcations of Vector Fields},
Springer, New York 1983.

\bibitem{HKW81}
B.D. Hassard, N.D. Kazarinoff, and Y.-H. Wan.
{\em Theory and Applications of Hopf Bifurcation},
London Math. Soc. Lecture Notes {\bf 41},
Cambridge University Press, Cambridge 1981.

\bibitem{HLC13}
X. He, W. Lu, and T. Chen.
On transverse stability of random dynamical system,
{\em Discrete and Continuous Dyn. Sys.} {\bf 33} (2013) 701-721.

\bibitem{HPS77}
M.W.~Hirsch, C.C.~Pugh, and M.~Shub.
{\em Invariant Manifolds}, Lect. Notes in Math. {\bf 583},
Springer, New York 1977.

\bibitem{HS74}
M.W. Hirsch and S. Smale.
{\em Differential Equations, Dynamical Systems, and Linear Algebra},
Academic Press, New York 1974.

\bibitem{IB18}
E.J. Izquierdo and R.D. Beer.  
From head to tail: a neuromechanical model
of forward locomotion in {\em Caenorhabditis elegans}, 
{\em Phil. Trans. R. Soc. Lond.} B {\bf 373} (2018); doi: 10.1098/rstb.2017.0374.

\bibitem{J62}
N. Jacobson.
{\em Lie Algebras}, Wiley, New York 1962.
 
\bibitem{JC19}
M. Jha and N.R. Chauhan.
A review on snake-like continuum robots for medical surgeries,
{\em IOP Conf. Ser.: Mater. Sci. Eng.} {\bf 691} (2019) 012093.

\bibitem{KH95}
A. Katok and B. Hasselblatt.
{\em Introduction to the Modern Theory of Dynamical Systems},
Cambridge University Press, Cambridge 1995.

\bibitem{K67}
A. Kelley. The stable, center-stable, center, center-unstable, unstable manifolds,
{\em J. Diff. Eq.} {\bf 3} (1967) 546--570.

\bibitem{K88}
N. Kopell.
Towards a theory of modelling central pattern generators,
in: {\em Neural Control of Rhythmic Movements in Vertebrates} (eds.
A.H. Cohen, S. Rossignol, and S. Grillner), Wiley, New York 1988.

\bibitem{KE88}
N. Kopell and G.B. Ermentrout.
Coupled oscillators and the design of central pattern generators,
{\em Math. Biosci.} {\bf 89} (1988) 14--23.

\bibitem{KE90}
N. Kopell and G.B. Ermentrout.
Phase transitions and other phenomena in chains of oscillators,
{\em SIAM J. Appl. Math.} {\bf 50} (1990) 1014--1052.

\bibitem{KF99}
W.A. Kunze and J.B. Furness.
The enteric nervous system and regulation of intestinal motility,
{\em Ann. Rev. Physiol.} {\bf 61} (1999) 117--42. 

\bibitem{K63}
I. Kupka. Contribution \`a la th\'eorie des champs g\'en\'eriques,
{\em Contrib. Diff. Eqs.} {\bf 2} (1963) 457--484; 
{\bf 3} (1964) 411--420.

\bibitem{K84}
Y.~Kuramoto. {\em Chemical Oscillations, Waves, and Turbulence},
Springer, Berlin 1984.

\bibitem{LC02}
C.R. Laing and C.C. Chow.
A spiking neuron model for binocular rivalry,
{\em J. Comput. Neurosci.} {\bf 12} (2002) 39--53.

\bibitem{L64}
J.P. LaSalle.
Recent advances in Liapunov stability theory, {\em SIAM Review} {\bf 6}
(1964) 1--11.

\bibitem{LL61}
J.P. LaSalle and S. Lefschetz. {\em Stability by Lyapunov's Second Method with Applications},
Academic Press, New York 1961.
 
\bibitem{L92}
A.M. Liapunov. 
{\em Obshchaya Zadacha Ustoichivosti Drizheniya}, Kharkov 1892,
and {\em Comm. Soc. Math. Kharkov} {\bf 3} (1893) 265--272;
French translation: Probl\`eme g\'en\'erale de la stabilit\'e du mouvement,
{\em Ann. Fac. Sci. Toulouse} {\bf 9} (1907) 203--474;
reproduced as Annals of Mathematics Studies {\bf 17},
Princeton University Press, Princeton 1947; English translation:
{\em Stability of Motion}, Academic Press, New York 1966.

\bibitem{L03}
D. Liberzon.
{\em Switching in Systems and Control}, Springer, New York 2003.

\bibitem{MY60}
L. Markus and H. Yamabe. Global stability criteria for differential
systems, {\em Osaka J. Math.} {\bf 12} (1960) 305--317.

\bibitem{MHKDA03}
C. Mehring, U. Hehl, M. Kubo, M. Diesmann, and A. Aertsen.
Activity dynamics and propagation of synchronous spiking in locally
connected random networks,
{\em Biol. Cybern.} {\bf 88} (2003) 395--408 .

\bibitem{MBBP19}
T. Menara, G. Baggio, D.S. Bassett, and F. Pasqualetti.
Stability conditions for cluster synchronization in
networks of heterogeneous Kuramoto oscillators,
arXiv:1806.06083v2 (2019).

\bibitem{M85}
J. Milnor. On the concept of attractor,
{\em Commun. Math. Phys.}{\bf 99} (1985) 177--195.

\bibitem{MLM20}
F. Morone, I. Leifer, and H.A. Makse.
Fibration symmetries uncover the building blocks of biological networks,
{\em Proc. Nat. Acad. Sci.} {\bf 117} (2020) 8306--  8314.

\bibitem{MWXLW17}
Z. Mu, H. Wang, W. Xu, T. Liu, and H. Wang.
Two types of snake-like robots for complex environment exploration: Design, development, and experiment,
{\em Adv. Mech. Eng.}{\bf 9} (2017);
doi: 10.1177/1687814017721.

\bibitem{MLS93}
R.M. Murray, Z. Li, and S.S. Sastry.
{\em A Mathematical Introduction to Robotic Manipulation},
CRC Press, Boca Raton 1993.

\bibitem{OIB21}	
E. Olivares, E.J. Izquierdo, and R.D. Beer.
A neuromechanical model of multiple network rhythmic pattern
generators for forward locomotion
in {\em C. elegans}, {\em Front. Comput. Neurosci.} {\bf 18} (2021); doi: 10.3380/fncom.2021.572339.

\bibitem{PCKHS92}
U. Parlitz, L.O. Chua, L. Kocarev, K.S. Halle, and A. Shang.
Transmission of digital signals by chaotic synchronization, 
{\em Int. J. Bif. Chaos} {\bf 2} (1992) 973--977. 

\bibitem{PC90}
L.M. Pecora and T.L. Carroll.
Synchronization in chaotic systems, {\em Phys. rev.Lett.} {\bf 64} (1990) 
821--825.

\bibitem{PC98}
L.M.~Pecora and T.L.~Carroll. Master stability functions for synchronized
coupled systems, {\em Phys. Rev. Lett.} {\bf 80} (1998) 2109--2112.

\bibitem{PCJM97}
L.M. Pecora, T.L. Carroll, G.A. Johnson, and D.J. Mar.
Fundamentals of synchronization in chaotic systems, concepts, and applications, 
{\em Chaos} {\bf 7} (1997) 520; doi: 10.1063/1.166278.

\bibitem{PC95}
G. P\'erez and H. A. Cerderia.
Extracting messages masked by chaos,
{\em Phys. Rev. Lett.} {\bf 74} (1995) 1970--.

\bibitem{PRK01}
A. Pikovsky, M. Rosenblum, and J. Kurths. 
{\em Synchronization: A Universal Concept in Nonlinear Sciences}, 
Cambridge University Press, Cambridge 2001.

\bibitem{PG06}
C.A. Pinto and M. Golubitsky. 
Central pattern generators for bipedal locomotion, 
{\em J. Math. Biol.} {\bf 53} (2006) 474--489.

\bibitem{PST93}
N. Platt, E. A. Spiegel, and C. Tresser.
On-off intermittency: A mechanism for bursting,
{\em Phys. Rev. Lett.} {\bf 70} (1993) 279--282.

\bibitem{PK17}
B.T. Polyak and Ya.I. Kvinto.
Stability and synchronization of oscillators:
new Lyapunov functions, 
{\em Automation and Remote Control} {\bf 78} (2017) 1234--1242;
original Russian text {\em Avtomatika i Telemekhanika} {\bf 7} (2017) 76--85.

\bibitem{SSSIT21}
K. Sakamoto, Z. Soh, M. Suzuki, Y. Iino, and T. Tsuji. 
Forward and backward locomotion patterns in {\em  C. elegans} generated
by a connectome‑based model simulation,
{\em Nature Scientific Reports} {\bf 11} (2021) 13737; doi: 10.1038/s41598-021-92690-2.

\bibitem{S02}
B.S.W.~Schr\"oder.
{\em Ordered Sets: An Introduction}, Birkh\"auser, Boston 2002.

\bibitem{SOWRK10}
S. Seok, C. D. Onal, R. Wood, D. Rus, and S. Kim.
Peristaltic locomotion with antagonistic actuators in soft robotics.
{\em 2010 IEEE International Conference on Robotics and Automation} (2010)
1228--1233; doi: 10.1109/ROBOT.2010.5509542.

\bibitem{SMG18}
H. Setareh, M. Deger, and W. Gerstner.
Excitable neuronal assemblies with adaptation as a building block of brain circuits
for velocity-controlled signal propagation,
{\em PLOS Comput. Biol.} {\bf 14} (2018) e1006216; doi: 0.1371/journal.pcbi.1006216.

\bibitem{SCRR07}
A. Shpiro, R. Curtu, J. Rinzel and N. Rubin.
Dynamical characteristics common to neuronal competition models,
{\em  J Neurophysiol} {\bf 97} (2007) 462--473.

\bibitem{S63}
S. Smale.
Stable manifolds for differential equations and diffeomorphisms,
{\em Ann. Scuola Normale Superiore Pisa} {\bf 17} (1963) 97--116.

\bibitem{S67}
S. Smale.
Differentiable dynamical systems,
{\em Bull. Amer. Math. Soc.} {\bf 73} (1967) 747--817. 

\bibitem{S14}
I. Stewart.
Symmetry-breaking in a rate model for a biped locomotion central pattern generator,
{\em Symmetry} {\bf 6} (2014) 23--66.

\bibitem{S20overdet}
I. Stewart.
Overdetermined ODEs and rigid periodic states in network dynamics, {\em Portugaliae Mathematica}, to appear; arxiv.org/abs/2112.15415 (2022).

\bibitem{SGP03}
I. Stewart, M. Golubitsky, and M. Pivato. Symmetry groupoids 
and patterns of synchrony in coupled cell networks, 
{\em SIAM J. Appl. Dynam. Sys.} {\bf 2} (2003) 609--646. 

\bibitem{SP08}
I. Stewart and M. Parker.
Periodic dynamics of coupled cell networks II: cyclic symmetry, {\em Dynamical Systems}
{\bf 23} (2008) 17--41.

\bibitem{SW23} 
I. Stewart and D. Wood. 
Stable synchronous propagation of signals
by feedforward networks: Examples and Applications, in preparation 2023.

\bibitem{TBB00}
E.A. Thomas, P.P. Bertrand, and J.C. Bornstein.  A computer
simulation of recurrent, excitatory networks of
sensory neurons of the gut in guinea-pig. {\em Neurosci.
Lett.} {\bf 287} (2000) 137--140.

\bibitem{WC72}
H.R.~Wilson and J.D.~Cowan.
Excitatory and inhibitory interactions in localized populations of model neurons,
{\em Biophys. J.} {\bf 12} (1972) 1--24.

\bibitem{Y02}
L.-S. Young.
What are SRB measures, and which dynamical systems have them?,
{\em J. Stat. Phys.} {\bf 108} (2002) 733--754.

\bibitem{ZHX20}
 Y. Zhong, L. Hu, and Y. Xu.
 Recent advances in design and actuation of continuum robots for medical applications,
{\em Actuators} {\bf 9} (2020) 142; doi: 10.3390/act9040142. 

\end{thebibliography}

\end{document}